\documentclass[10pt,english,american]{article}
\usepackage{ae,aecompl}
\usepackage{srcltx}

\usepackage[T1]{fontenc}
\usepackage[latin9]{inputenc}
\usepackage{placeins}
\usepackage{geometry}
\usepackage{color}
\usepackage{float}
\usepackage{units}
\usepackage{bm}
\usepackage{amsmath}
\usepackage{amssymb}
\usepackage{stmaryrd}
\usepackage{graphicx}
\usepackage{setspace}
\makeatletter

\providecommand{\tabularnewline}{\\}
\floatstyle{ruled}
\newfloat{algorithm}{tbp}{loa}
\providecommand{\algorithmname}{Algorithm}
\floatname{algorithm}{\protect\algorithmname}

\newcommand{\lyxaddress}[1]{
	\par {\raggedright #1
	\vspace{1.4em}
	\noindent\par}
}

\@ifundefined{date}{}{\date{}}
\usepackage{bm,url}
\providecommand{\dif}{\mathrm{d}}
\usepackage{graphicx,wasysym,datetime2}
\usepackage[T1]{fontenc}
\usepackage{calrsfs}

\@ifundefined{showcaptionsetup}{}{%
 \PassOptionsToPackage{caption=false}{subfig}}
\usepackage{subfig}
\makeatother

\usepackage{babel}
\addto\captionsenglish{\renewcommand{\algorithmname}{Algorithm}}

\begin{document}
\title{\textbf{\emph{Continuous gap contact formulation based on the
    screened Poisson equation}}}

\author{\selectlanguage{english}%
P. Areias$^{a,b\star}$, \textcolor{black}{N. Sukumar}$^{c}$ and J.
Ambr\'osio$^{a,b}$}
\maketitle
\selectlanguage{english}%

\lyxaddress{\begin{center}
{\footnotesize{}$^{a}$}\textbf{\emph{\footnotesize{}DEM}}{\footnotesize{}
- Departamento de Engenharia Mec\^anica }\\
\textbf{\footnotesize{}Instituto Superior T�cnico}{\footnotesize{}}\\
{\footnotesize{}Universidade de Lisboa}\\
{\footnotesize{}Building Mechanics II, Avenida Rovisco Pais 1, 1049-001}\\
{\footnotesize{} Lisboa, Portugal\vspace{.3cm}}\\
$^{b}$\textbf{\emph{\footnotesize{}IDMEC}}{\footnotesize{} - Instituto
Superior T�cnico}\\
{\footnotesize{} Avenida Rovisco Pais 1, 1049-001}\\
{\footnotesize{}Lisboa, Portugal\vspace{.3cm}}\\
{\footnotesize{}$^{c}$}\textbf{\emph{\footnotesize{}UC-Davis }}{\footnotesize{}-}
\textbf{\footnotesize{}University of California at Davis}{\footnotesize{}}\\
{\footnotesize{}Department of Civil \& Environmental Engineering}\\
{\footnotesize{}One Shields Avenue}\\
{\footnotesize{}Davis, CA 95616. U.S.A}\vspace{.3cm}\\
{\footnotesize{}$^{\star}$pedro.areias@tecnico.ulisboa.pt}
\par\end{center}}

\begin{abstract}
 We introduce a PDE-based node-to-element contact formulation as an
alternative to classical, purely geometrical formulations. It is challenging
to devise solutions to nonsmooth contact problem with continuous gap
 using finite element discretizations. 
We herein achieve this objective by constructing  an
approximate distance function (ADF) to the boundaries of solid objects, and in
doing so, also obtain universal uniqueness
of contact detection. Unilateral constraints are implemented using
a mixed model combining the screened Poisson equation and a force
element, which has the topology of a continuum element containing
an additional incident node. An ADF is obtained by solving the screened Poisson equation with constant essential
boundary conditions and  a variable transformation. The ADF does not
explicitly depend on the number of objects and a single solution of
the partial differential equation for this field uniquely defines
the contact conditions for all incident points in the mesh. Having
an ADF field to any obstacle circumvents the multiple target surfaces and 
avoids the specific data structures present in traditional contact-impact 
algorithms.  We also relax the interpretation of the Lagrange multipliers
as contact forces, and the Courant--Beltrami function is used with a mixed
formulation producing the required differentiable result. We demonstrate
the advantages of the new approach in two- and three-dimensional problems
that are solved using Newton iterations. Simultaneous constraints for each incident point are considered.
\end{abstract}


\maketitle

\section{Introduction}
Considerable interest has recently emerged with the enforcement of
essential boundary conditions~\cite{sukumar2022} and contact detection
algorithms~\cite{lai2022} using approximate distance functions (ADFs).
Realistic contact frameworks for finite-strain problems make use of
three classes of algorithms that are often coupled (see~\cite{wriggerscontact}):
\begin{enumerate}
\item Techniques of unilateral constraint enforcement. Examples of these
are direct elimination (also reduced-gradient algorithms), penalty
and barrier methods, and Lagrange multiplier methods with the complementarity
condition~\cite{simo1992c}. Stresses from the underlying discretization
are often used to assist the normal condition with Nitsche's method~\cite{nitsche1970}.
\item Frictional effects (and more generally,
constitutive contact laws) 
and cone-complementarity forms~\cite{kanno2006,areias2015c}.
Solution paradigms are augmented Lagrangian methods for friction~\cite{laursen1993},
cone-projection techniques~\cite{desaxce1998} and the so-called yield-limited
algorithm for friction~\cite{jones2000}.
\item Discretization methods. Of concern are distance calculations, estimation
of tangent velocities and general discretization arrangements. In
contemporaneous use are (see~\cite{wriggerscontact}) node-to-edge/face,
edge/face-to-/edge/face and mortar discretizations.
\end{enumerate}

In terms of contact enforcement and friction algorithms, finite displacement
contact problems are typically addressed with well-established contact
algorithms, often derived from solutions developed by groups working
on variational inequalities and nonsmooth analysis. However, in the
area of contact discretization and the related area of contact kinematics
\cite{simo1985}, there are still challenges to be addressed in terms
of ensuring the complete robustness of implicit codes. One of the
earliest papers on contact discretization was by Chan and Tuba \cite{chan1971},
who considered contact with a plane and used symmetry to analyze cylinder-to-cylinder
contact. Francavilla and Zienkiewicz \cite{francavilla1975} later
proposed an extension to node-to-node contact in small strains, allowing
for direct comparison with closed-form solutions. The extension to
finite strains requires further development, and the logical development
was the node-to-segment approach, as described in the work of Hallquist
\cite{hallquist1985}. Although node-to-segment algorithms are known
to entail many defficiencies, most of the drawbacks have been addressed.
Jumps and discontinuities resulting from nodes sliding between edges
can be removed by smoothing and regularization \cite{neto2014}. Satisfaction
of patch-test, which is fundamental for convergence, can be enforced
by careful penalty weighing \cite{zavarise2009b,zavarise2009}. It is well known
that single-pass versions of the node-to-segment method result in
interference and double-pass can produce mesh interlocking, see \cite{puso2004,puso2004b}.
This shortcoming has eluded any attempts of a solution and has motivated
the development of surface-to-surface algorithms. One of the first
surface-to-surface algorithms was introduced by Simo, Wriggers, and
Taylor \cite{simo1985}. Zavarise and Wriggers \cite{zavarise1998}
presented a complete treatment of the 2D case and further developments
were supported by parallel work in nonconforming meshes, see \cite{kim2001}.
A review article on mortar methods for contact problems \cite{laursen2012}
where stabilization is discussed and an exhaustive detailing of most
geometric possibilities of contact was presented by Farah, Wall and
Popp \cite{farah2018}. This paper reveals the significant effort
that is required to obtain a robust contact algorithm. An interesting
alternative approach to contact discretization has been proposed by
Wriggers, Schr\"oder, and Schwarz \cite{wriggers2013}, who use an intermediate
mesh with a specialized anisotropic hyperelastic law to represent
contact interactions. In the context of large, explicit codes, Kane
et al. \cite{kane1999} introduced an interference function, or gap,
based on volumes of overlapping, allowing non-smooth contact to be
dealt in traditional smooth-based algorithms. 

In addition to these general developments, there have been specialized
algorithms for rods, beams, and other structures. Litewka et al. \cite{litewka2013}
as well as Neto et al. \cite{neto2016,neto2017}, have produced efficient
algorithms for beam contact. For large-scale analysis of beams, cables
and ropes, Meier et al. \cite{meier2017} have shown how beam specialization
can be very efficient when a large number of entities is involved.

Recently, interest has emerged on using approximate distance 
functions~\cite{wolff2013,liu2020cm,aguirre2020,lai2022} as alternatives to algorithms 
that heavily rely on computational geometry. These algorithms
resolve the non-uniqueness of the projection, but still require 
geometric computations. 
In~\cite{wolff2013}, Wolff and Bucher proposed a
local construction to obtain distances inside any object, but still
require geometric calculations, especially for the integration scheme
along surfaces. Liu et al.~\cite{liu2020cm} have combined finite
elements with distance potential discrete element method (DEM) in
2D within an explicit integration framework. A geometric-based distance
function is constructed and contact forces stem from this construction.
In~\cite{aguirre2020}, the analysis of thin rods is performed using
classical contact but closed-form contact detection is achieved by
a signed-distance function defined on a voxel-type grid. In~\cite{lai2022},
a pre-established set of shapes is considered, and a
function is defined for each particle in a DEM (discrete element
method) setting with a projection that follows. 
In the context
of computer graphics and computational geometry, Macklin et al.~\cite{macklin2020}
introduced an element-wise local optimization algorithm to determine
the closest-distance between the signed-distance-function 
isosurface and face elements. Although close to what is proposed here,
no solution to a partial differential equation (PDE) 
is proposed and significant geometric calculations are still required.

In this paper, we adopt a different method, which aims to be more general
and less geometric-based. This is possible due to the equivalence
between the solution of the Eikonal equation and the distance 
function~\cite{belyaev2015}. 
It is worth noting that very efficient linear algorithms are available
to solve regularized Eikonal equations, as discussed by Crane 
et al.~\cite{crane2013}. The work in~\cite{crane2013}
provides a viable solution 
for contact detection in computational mechanics. Applications 
extend to beyond contact mechanics and they provide a solution for the 
long-standing issue of imposing essential boundary conditions in 
meshfree methods~\cite{sukumar2022}.

We solve a partial differential equation (PDE) that produces an ADF
(approximate distance function) that converges to the distance function
as a length parameter tends to zero. The relation between the screened
Poisson equation (also identified as Helmholtz-like equation), which
is adopted in computational damage and fracture~\cite{peerlings,peerlings2001}
and the Eikonal equation is well understood~\cite{guler2014}. Regularizations
of the Eikonal equation such as the geodesics-in-heat~\cite{crane2013} and the
screened Poisson equation are discussed by Belyaev and Fayolle~\cite{belyaev2020}.
We take advantage of the latter for computational contact mechanics.
Specifically, the proposed algorithm solves well-known shortcomings
of geometric-based contact enforcement:
\begin{enumerate}
\item Geometric calculations are reduced to the detection of a target element
for each incident node.
\item The gap function $g(\bm{x})$ is continuous, since the solution of
the screened Poisson equation is a continuous function.
\item The contact force direction is unique and obtained as the gradient
of $g(\bm{x})$. 
\item Since the Courant--Beltrami penalty function is adopted, contact force
is continuous in the normal direction.
\end{enumerate}
Concerning the importance of solving the uniqueness problem, Konyukhov
and Schweizerhof~\cite{konyukhov2008} have shown that considerable
computational geometry must be in place to ensure uniqueness and existence
of node-to-line projection. Another important computational effort
was presented by Farah et al.~\cite{farah2018} to geometrically address
all cases in 3D with mortar interpolation. Extensions to higher dimensions
require even more intricate housekeeping. When compared with the geometrical
approach to the distance function~\cite{aguirre2020,lai2022}, the
algorithm is much simpler at the cost of a solution of an additional
PDE. Distance functions can be generated by level set solutions of
the transport equation~\cite{russo2000}.

The remainder of this paper is organized as follows. In Section \ref{sec:Approximate-distance-function},
the approximate distance function is introduced as the solution of
a regularization of the Eikonal equation. In Section \ref{sec:Discretization},
details concerning the discretization are introduced. The overall
algorithm, including the important step-control, is presented in Section
\ref{sec:Algorithm-based-on}. Verification and validation examples
are shown in Section \ref{sec:Numerical-tests} in both 2D and 3D.
Finally, in Section \ref{sec:Conclusions}, we present the advantages
and shortcomings of the present algorithm, as well as suggestions
for further improvements.

\section{Approximate distance function (ADF)}

\label{sec:Approximate-distance-function}

Let $\Omega\subset\mathbb{R}^{d}$ be a deformed configuration of
a given body in $d$-dimensional space and $\Omega_{0}$ the corresponding undeformed
configuration. The boundaries of these configurations are $\Gamma$
and $\Gamma_{0}$, respectively. Let us consider deformed coordinates
of an arbitrary point $\text{\ensuremath{\bm{x}\in\Omega}}$ and a
specific point, called incident, with coordinates $\bm{x}_{I}$. For
reasons that will become clear, we also consider a \emph{potential
}function $\phi\left(\bm{x}_{I}\right)$,  which is the solution of
a scalar PDE.  We are concerned with an approximation of the signed-distance
function. The so-called \emph{gap function }is now introduced as a
differentiable function $g\left[\phi\left(\bm{x}_{I}\right)\right]$
such that~\cite{wriggerscontact}:
\begin{equation}
g\left[\phi\left(\bm{x}_{I}\right)\right]\begin{cases}
<0 & \quad\bm{x}_{I}\in\Omega\\
=0 & \quad\bm{x}_{I}\in\Gamma\\
>0 & \quad\bm{x}_{I}\notin\Omega\cup\Gamma
\end{cases}.
\label{eq:gtriple}
\end{equation}

If a unique normal $\bm{n}\left(\bm{x}_{I}\right)$ exists for $\bm{x}_{I}\in\Gamma$,
we can decompose the gradient of $g\left[\phi\left(\bm{x}_{I}\right)\right]$
into parallel ($\parallel$) and orthogonal ($\perp$) terms: $\nabla g\left[\phi\left(\bm{x}_{I}\right)\right]=\left\{ \nabla g\left[\phi\left(\bm{x}_{I}\right)\right]\right\} _{\parallel}+\left\{ \nabla g\left[\phi\left(\bm{x}_{I}\right)\right]\right\} _{\perp}$,
with $\left\{ \nabla g\left[\phi\left(\bm{x}_{I}\right)\right]\right\} _{\perp}\sslash\bm{n}\left(\bm{x}_{I}\right)$.
Since $g\left[\phi\left(\bm{x}_{I}\right)\right]=0$ for $\bm{x}_{I}\in\Gamma$ we
have $\nabla g\left[\phi\left(\bm{x}_{I}\right)\right]\sslash\bm{n}\left(\bm{x}_{I}\right)$
on those points. In the frictionless case, the normal contact force component 
 is identified as $f_{c}$ and contact conditions correspond to the
following complementarity conditions~\cite{simo1992c}:
\begin{align}
g\left[\phi\left(\bm{x}_{I}\right)\right]f_{c} & =0 , \nonumber \\
f_{c} & \geq 0 , \label{eq:gplu}\\
g\left[\phi\left(\bm{x}_{I}\right)\right] & \geq 0 , \nonumber 
\end{align}
or equivalently, $\left\langle gg\left[\phi\left(\bm{x}_{I}\right)\right]+f_{c}\right\rangle -f_{c}=0$
where $\left\langle x\right\rangle =\nicefrac{\left(x+|x|\right)}{2}$.
The vector form of the contact force is given by $\bm{f}_{c}=f_{c}\nabla g\left[\phi\left(\bm{x}_{I}\right)\right]$.
Assuming that a function $g\left[\phi\left(\bm{x}_{I}\right)\right]$, and its first
and second derivatives are available from the solution of a PDE, we
approximately satisfy conditions (\ref{eq:gplu}) by defining the
contact force as follows:
\begin{equation}
\bm{f}_{c}=f_{c}\nabla g\left[\phi\left(\bm{x}_{I}\right)\right]=\kappa\min\left\{ 0,g\left[\phi\left(\bm{x}_{I}\right)\right]\right\} ^{2}\nabla g\left[\phi\left(\bm{x}_{I}\right)\right]\qquad(\kappa>0), \label{eq:fcx}
\end{equation}
where $\kappa$ is a penalty parameter for the Courant--Beltrami
function~\cite{courant1943,beltrami1969}. The Jacobian of this force is given by:
\begin{equation}
\bm{J}=\kappa\min\left\{ 0,g\left[\phi\left(\bm{x}_{I}\right)\right]\right\} ^{2}\nabla
\otimes \nabla g\left[\phi\left(\bm{x}_{I}\right)\right]+2\kappa\min\left\{ 0,g\left[\phi\left(\bm{x}_{I}\right)\right]\right\} \nabla g\left[\phi\left(\bm{x}_{I}\right)\right]\otimes\nabla g\left[\phi\left(\bm{x}_{I}\right)\right]. 
\label{eq:stiffn}
\end{equation}

Varadhan~\cite{varadhan1967} established that the solution of the
screened Poisson equation:
\begin{subequations}\label{eq:cln2}
\begin{align}
c_{L} \nabla^2 \phi(\bm{x}) - \phi(\bm{x}) &= 0\quad\textrm{in } \Omega \label{eq:cln2-a}\\
\phi (\bm{x}) & =1 \quad \textrm{on } \Gamma  \label{eq:cln2-b}
\end{align}
\end{subequations}
produces an ADF given by $-c_{L}\log\left[\phi\left(\bm{x}\right)\right]$. This property has been recently studied by Guler et al.~\cite{guler2014}
and Belyaev et al.~\cite{belyaev2015,belyaev2020}. The exact distance
is obtained as a limit:
\begin{equation}
d\left(\bm{x}\right)=-\lim_{c_{L}\rightarrow0}c_{L}\log\left[\phi\left(\bm{x}\right)\right].\label{eq:ddeterm}
\end{equation}
which is the solution of the Eikonal equation. Proof of this limit is provided in Varadhan~\cite{varadhan1967}.
We transform the ADF to a signed-ADF, by introducing the sign of outer
($+$) or inner ($-$) consisting of 
\begin{equation}
g\left(\bm{x}\right)=\pm c_{L}\log\left[\phi\left(\bm{x}\right)\right].\label{eq:gx}
\end{equation}
Note that if the plus sign is adopted in (\ref{eq:gx}), then inner
points will result in a negative gap $g(\bm{x}$). The gradient of
$g(\bm{x})$ results in 

\begin{equation}
\nabla g\left(\bm{x}\right)=\pm\frac{c_{L}}{\phi\left(\bm{x}\right)}\nabla\phi\left(\bm{x}\right)\label{eq:nablag}
\end{equation}
\\
and the Hessian of $g(\bm{x})$ is obtained in a straightforward manner:
\begin{equation}
\nabla\otimes\nabla\left[g\left(\bm{x}\right)\right]=\pm\frac{c_{L}}{\phi\left(\bm{x}\right)}\left\{ \nabla\otimes\nabla\left[\phi\left(\bm{x}\right)\right]-\frac{1}{\phi\left(\bm{x}\right)}\nabla\phi\left(\bm{x}\right)\otimes\nabla\phi\left(\bm{x}\right)\right\}.
 \label{eq:n2g}
\end{equation}
Note that $c_L$ is a square of a characteristic length, i.e. $c_L=l_c^2$, which is here taken as a solution parameter.

Using a test field $\delta\phi(\bm{x})$, the weak form
of (\ref{eq:cln2}) is written as 
\begin{equation}
\int_{\Omega}c_{L}\nabla\phi\left(\bm{x}\right)\cdot\nabla\delta\phi\left(\bm{x}\right)\dif V+\int_{\Omega}\phi\left(\bm{x}\right)\delta\phi\left(\bm{x}\right)\dif V=\int_{\Gamma}c_{L}\delta\phi\,\nabla\phi\left(\bm{x}\right)\cdot\bm{n}\left(\bm{x}\right)\dif A , 
\label{eq:weakclphi}
\end{equation}
where $\textrm{d}A$ and $\textrm{d}V$ are differential area and volume
elements, respectively.
Since an essential boundary condition is imposed on
$\Gamma$ such that $\phi\left(\bm{x}\right)=1$ for $x\in\Gamma$,
it follows that 
$\delta \phi(\bm{x}) = 0$ on $\Gamma$ and the right-hand-side 
of (\ref{eq:weakclphi}) vanishes. 

\section{Discretization}

\label{sec:Discretization}

In an interference condition, each interfering node, with coordinates
$\bm{x}_{N}$, will fall within a given continuum element. The parent-domain
coordinates $\bm{\xi}$ for the incident node $\bm{x}_{I}$ also depends
on element nodal coordinates. Parent-domain coordinates are given
by:
\begin{subequations}
\begin{align}
\bm{\xi}_{I} &= \underset{\bm{\xi}}{\textrm{argmin}}\left[\left\Vert \bm{x}\left(\bm{\xi}\right)-\bm{x}_{I}\right\Vert \right], \label{eq:xi} \\
\intertext{and it is straightforward to show that for a triangle,}
\xi_{I1} &= \frac{x_{I2}\left(x_{11}-x_{31}\right)+x_{12}x_{31}-x_{11}x_{32}+x_{I1}\left(x_{32}-x_{12}\right)}{x_{11}\left(x_{22}-x_{32}\right)+x_{12}\left(x_{31}-x_{21}\right)+x_{21}x_{32}-x_{22}x_{31}}, \label{eq:xi1} \\
\xi_{I2} &=\frac{x_{I2}\left(x_{21}-x_{11}\right)+x_{12}x_{21}-x_{11}x_{22}+x_{I1}\left(x_{12}-x_{22}\right)}{x_{11}\left(x_{22}-x_{32}\right)+x_{12}\left(x_{31}-x_{21}\right)+x_{21}x_{32}-x_{22}x_{31}}, \label{eq:xi2}
\end{align}
\end{subequations}
with similar expressions for a tetrahedron~\cite{Areias2022cgithub}.
The continuum element interpolation is as follows:
\begin{equation}
\bm{x}\left(\bm{\xi}\right)=\sum_{K=1}^{d+1}N_{K}\left(\bm{\xi}\right)\bm{x}_{K}=\bm{N}\left(\bm{\xi}\right)\cdot\bm{x}_{N},\label{eq:xxi}
\end{equation}
where $N_{K}\left(\bm{\xi}\right)$ with $K=1,\ldots,d+1$
are the shape functions of a triangular or tetrahedral element. 
Therefore, \eqref{eq:xi} can be written as:
\begin{equation}
\bm{\xi}_{I} = \underset{\bm{\xi}}{\textrm{argmin}}\left[\left\Vert \bm{N}\left(\bm{\xi}\right)\cdot\bm{x}_{N}-\bm{x}_{I}\right\Vert \right]=\bm{a}_N\left(\bm{x}_C\right)\label{eq:firstian}
\end{equation}
In (\ref{eq:firstian}), we group the continuum node coordinates and the
incident node coordinates in a single array $\bm{x}_{C}=\left\{ \begin{array}{cc}
\bm{x}_{N} & \bm{x}_{I}\end{array}\right\} ^{T}$with cardinality $\#\bm{x}_{C}=d\left(d+2\right)$.  We adopt the notation $\bm{x}_{N}$ for the coordinates of the continuum
element. For triangular and tetrahedral discretizations, $\bm{\xi}_{I}$
is a function of $\bm{x}_{N}$ and $\bm{x}_{I}$:

\begin{equation}
\bm{\xi}_{I}=\bm{a}_{N}\left(\left\{ \begin{array}{c}
\bm{x}_{N}\\
\bm{x}_{I}
\end{array}\right\} \right)=\bm{a}_{N}\left(\bm{x}_{C}\right). \label{eq:xian}
\end{equation}
The first and
second derivatives of $\bm{a}_{N}$ with respect to $\bm{x}_{C}$
make use of the following notation:
\begin{align}
\bm{A}_{N}\left(\bm{x}_{C}\right) & =\frac{\dif\bm{a}_{N}\left(\bm{x}_{C}\right)}{\dif\bm{x}_{C}}\qquad 
\left( d\times\left[d\left(d+2\right)\right] \right), \label{eq:firstderivative}\\
\mathcal{A}_{N}\left(\bm{x}_{C}\right) & =\frac{\dif\bm{A}_{N}\left(\bm{x}_{C}\right)}{\dif\bm{x}_{C}}\qquad 
\left( d\times\left[d\left(d+2\right)\right]^{2} \right). \label{eq:secondderivative}
\end{align}
Source code for these functions is available in \texttt{Github}~\cite{Areias2022cgithub}. A mixed formulation is adopted with equal-order interpolation for
the displacement $\bm{u}$ and the function $\phi$. For a set of
nodal displacements $\bm{u}_{N}$ and nodal potential values $\bm{\phi}_{N}$:
\begin{subequations}
\begin{align}
\bm{u}\left(\bm{x}_{C}\right) & =\bm{N}_{u}\left(\bm{x}_{C}\right)\cdot\bm{u}_{N},
\label{eq:uxi}\\
\phi\left(\bm{x}_{C}\right) & =\bm{N}_{\phi}\left(\bm{x}_{C}\right)\cdot\bm{\phi}_{N},
\label{eq:phixi} \\
\intertext{where in three dimensions}
\bm{N}_{u}\left(\bm{x}_{C}\right) & =\left[\begin{array}{ccccc}
\cdots & N_{K}\left[\bm{a}_{N}\left(\bm{x}_{C}\right)\right] & 0 & 0 & \cdots\\
\cdots & 0 & N_{K}\left[\bm{a}_{N}\left(\bm{x}_{C}\right)\right] & 0 & \cdots\\
\cdots & 0 & 0 & N_{K}\left[\bm{a}_{N}\left(\bm{x}_{C}\right)\right] & \cdots
\end{array}\right],\label{eq:shu}\\
\bm{N}_{\phi}\left(\bm{x}_{C}\right) & =\left[\begin{array}{ccc}
\cdots & N_{K}\left[\bm{a}_{N}\left(\bm{x}_{C}\right)\right] & \cdots\end{array}\right]^{T}
,
\end{align}
\end{subequations}
where $\bm{\xi}_{I}=\bm{a}_{N}\left(\bm{x}_{C}\right)$. First and
second derivatives of $N_{K}\left[\bm{a}_{N}\left(\bm{x}_{C}\right)\right]$
are determined from the chain rule:
\begin{align}
\frac{\dif N_{K}}{\dif\bm{x}_{C}} & =\frac{\dif N_{K}}{\dif\bm{\xi}_{I}}\cdot\bm{A}_{N}\left(\bm{x}_{C}\right),\label{eq:dndxc}\\
\frac{\dif^{2}N_{K}}{\dif\bm{x}_{C}^{2}} & =\bm{A}_{N}^{T}\left(\bm{x}_{C}\right)\cdot\frac{\dif^{2}N_{K}}{\dif\bm{\xi}_{I}^{2}}\cdot\bm{A}_{N}\left(\bm{x}_{C}\right)+\frac{\dif N_{K}}{\dif\bm{\xi}_{I}}\cdot\mathcal{A}_{N}\left(\bm{x}_{C}\right).
\end{align}

For the test function of the incident point, we have
\begin{equation}
\delta\phi\left(\bm{x}_{C}\right) =\bm{N}_{\phi}\left(\bm{x}_{C}\right)\cdot\delta\boldsymbol{\phi}_{N}+\bm{\phi}_{N}\cdot\frac{\dif\bm{N}_{\phi}\left(\bm{x}_{C}\right)}{\dif\bm{x}_{C}}\cdot\delta\bm{x}_{C}.\label{eq:dphi}
\end{equation}

For linear continuum elements, the
second variation of $\phi\left(\bm{\xi}\right)$ is given by the following rule:
\begin{align}
\dif\delta\phi\left(\bm{x}_{C}\right) =\, &\, \delta\bm{\phi}_{N}\cdot\frac{\dif\bm{N}_{\phi}\left(\bm{x}_{C}\right)}{\dif\bm{x}_{C}}\cdot\dif\bm{x}_{C}+\dif\bm{\phi}_{N}\cdot\frac{\dif\bm{N}_{\phi}\left(\bm{x}_{C}\right)}{\dif\bm{x}_{C}}\cdot\delta\bm{x}_{C} \nonumber \\
 & +\bm{\phi}_{N}\cdot\frac{\dif^{2}\bm{N}_{\phi}\left(\bm{x}_{C}\right)}{\dif\bm{x}_{C}^{2}}:\left(\delta\bm{x}_{C}\otimes\dif\bm{x}_{C}\right)\label{eq:sv}
\end{align}

Since the gradient of $\phi$ makes use of the continuum part of the
formulation, we obtain:

\begin{align}
\nabla\phi\left(\bm{\xi}\right) & =\bm{\phi}_{N}\cdot\underbrace{\frac{\dif\boldsymbol{N}_{\phi}\left(\bm{\xi}\right)}{\dif\bm{\xi}}\cdot\bm{j}^{-1}}_{\nicefrac{\dif\bm{N}_{\phi}}{\dif\bm{x}}},\qquad\nabla\delta\phi\left(\bm{\xi}\right)=\delta\bm{\phi}_{N}\cdot\frac{\dif\boldsymbol{N}_{\phi}\left(\bm{\xi}\right)}{\dif\bm{\xi}}\cdot\bm{j}^{-1}\label{eq:gphijac}
\end{align}
\\
where $\bm{j}$ is the Jacobian matrix in the deformed configuration.
The element residual and stiffness are obtained from these quantities
and available in \texttt{Github}~\cite{Areias2022cgithub}. Use is
made of the \texttt{Acegen}~\cite{korelc2002} add-on 
to \texttt{Mathematica}~\cite{mathematica}
to obtain the source code for the final expressions.

\section{Algorithm and step control}

\begin{figure}
\centering
\includegraphics[width=0.7\textwidth]{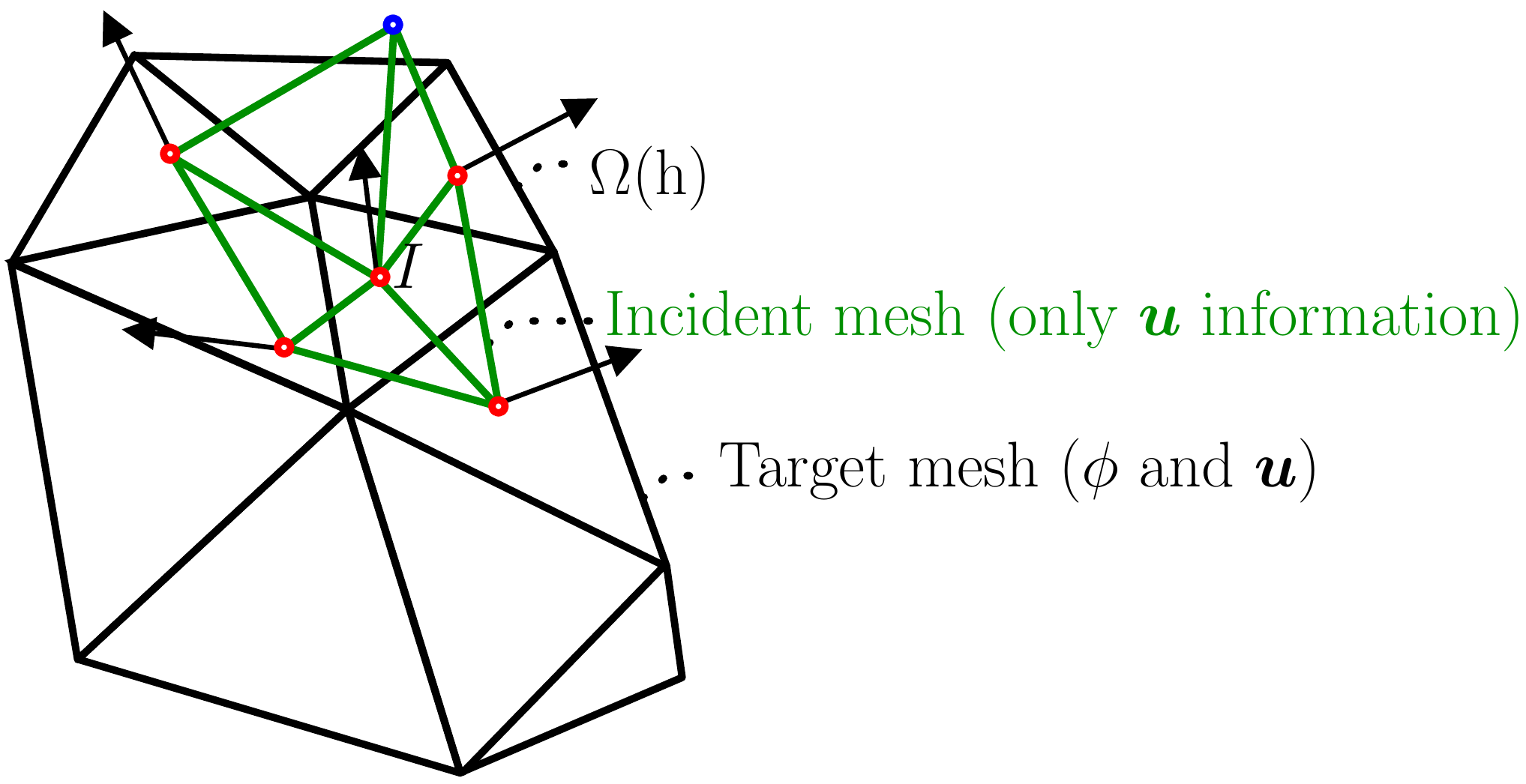}
\caption{Relevant quantities for the definition of a contact discretization 
         for element $e$.}\label{fig:Definition-of-a}
\end{figure}

All nodes are considered candidates and all elements are targets.
A simple bucket-sort strategy is adopted to reduce the computational
cost. In addition, 

Step control is required to avoid the change of target during
Newton-Raphson iteration. The screened Poisson equation is solved
for all bodies in the analyses. Figure~\ref{fig:Definition-of-a} shows
the simple mesh overlapping under consideration. The resulting algorithm
is straightforward. Main characteristics are:
\begin{itemize}
\item All nodes are candidate incident nodes.
\item All elements are generalized targets.
\end{itemize}
\label{sec:Algorithm-based-on}

The modifications required for a classical nonlinear Finite Element
software (in our case \texttt{SimPlas}~\cite{simplascode}) to include this
contact algorithm are modest. Algorithm~\ref{alg:Contact-algorithm-based}
presents the main tasks. In this Algorithm, blue boxes denote the
contact detection task, which here is limited to:
\begin{enumerate}
\item Detect nearest neighbor (in terms of distance) elements for each node.
A bucket ordering is adopted.
\item Find if the node is interior to any of the neighbor elements by use
of the shape functions for triangles and tetrahedra. This is performed
in closed form.
\item If changes occurred in the targets, update the connectivity table
and the sparse matrix addressing. Gustavson's algorithm~\cite{gustavson1978}
is adopted to perform the updating in the assembling process.
\end{enumerate}

In terms of detection, the following algorithm is adopted:
\begin{enumerate}
\item Find all exterior faces, as faces sharing only one tetrahedron.
\item Find all exterior nodes, as nodes belonging to exterior faces.
\item Insert all continuum elements and all exterior nodes in two bucket
lists. Deformed coordinates of nodes and deformed coordinates of element
centroids are considered.
\item Cycle through all exterior nodes
\begin{enumerate}
\item Determine the element bucket from the node coordinates
\item Cycle through all elements ($e$) in the $3^{3}=27$ buckets surrounding
the node
\begin{enumerate}
\item If the distance from the node to the element centroid is greater than
twice the edge size, go to the next element
\item Calculates the projection on the element ($\bm{\xi}_{I}$) and the
corresponding shape functions $\bm{N}\left(\bm{\xi}_{I}\right)$.
\item If $0\leq N_{K}\left(\bm{\xi}_{I}\right)\leq1$ then $e$ is the target
element. If the target element has changed, then flag the solver for connectivity update.
\end{enumerate}
\end{enumerate}
\end{enumerate}

Since the algorithm assumes a fixed connectivity table during Newton
iterations, a verification is required after each converged iteration
to check if targets have changed since last determined. If this occurs,
a new iteration sequence with revised targets is performed. 

\begin{algorithm}
\caption{\label{alg:Contact-algorithm-based}Staggered contact algorithm based on the
approximate distance function. $\lambda$ is the load/displacement factor.}

\begin{centering}
\includegraphics[width=0.6\textwidth]{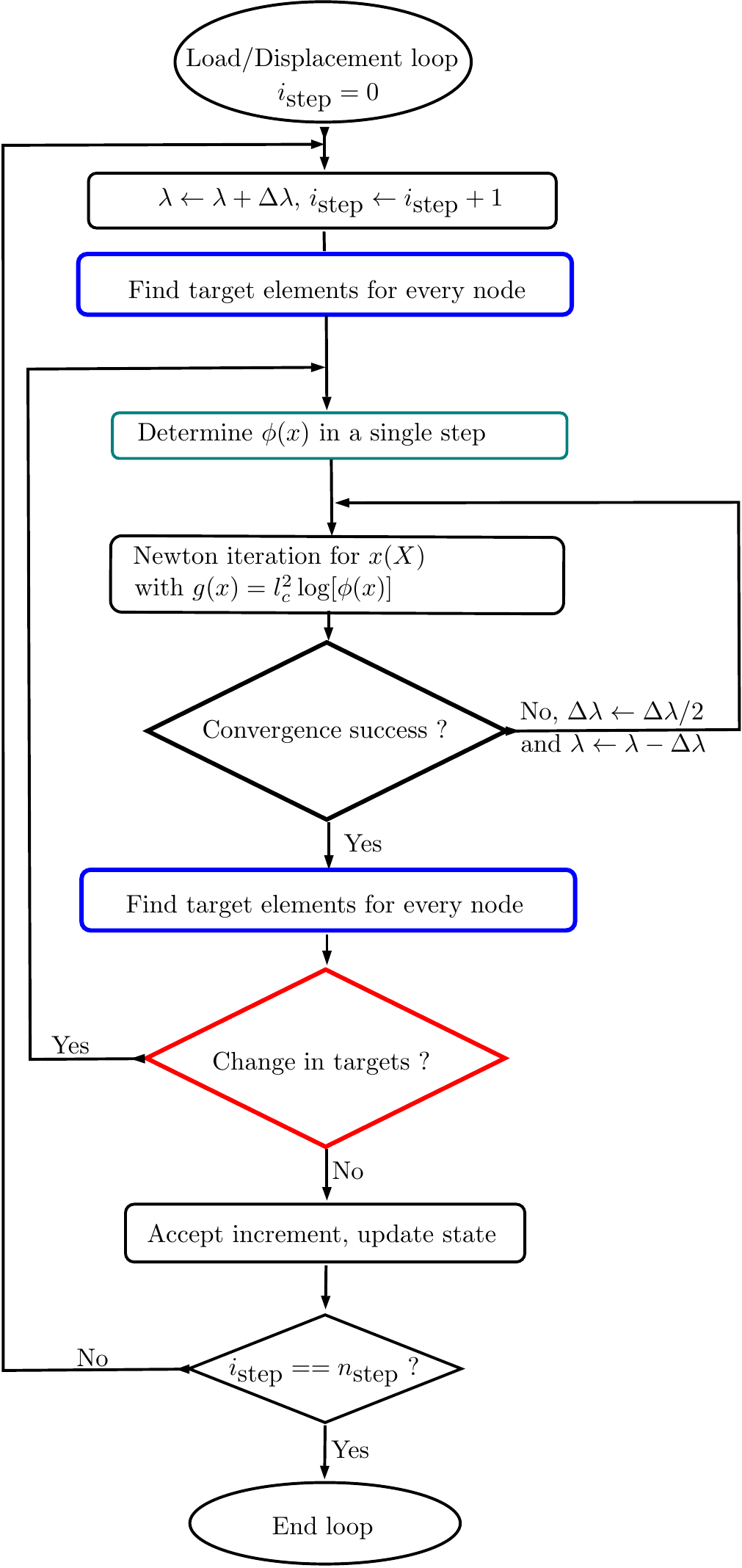}
\par\end{centering}
\end{algorithm}

The only modification required to a classical FEM code is the solution of the screened-Poisson equation, the green box in Algorithm \ref{alg:Contact-algorithm-based}. The cost of this solution is negligible when compared with the nonlinear solution process since the equation is linear and scalar. It is equivalent to the cost of a steady-state heat conduction solution. Note that this corresponds to a staggered algorithm.

\section{Numerical tests}
\label{sec:Numerical-tests}
Numerical examples are solved with our in-house software, 
\texttt{SimPlas}~\cite{simplascode}, using the new node-to-element contact element.
Only triangles and tetrahedra are assessed at this time, which provide
an exact solution for $\bm{\xi}_{I}$. \texttt{Mathematica}~\cite{mathematica}
with the add-on \texttt{Acegen}~\cite{korelc2002} is employed to obtain the
specific source code. All runs are quasi-static and make use of a Neo-Hookean model. If $\bm{C}$ is the right Cauchy-Green tensor, then

\[
\bm{S}=2\nicefrac{\dif\psi\left(\bm{C}\right)}{\dif\bm{C}}
\]
 where $\psi\left(\bm{C}\right)=\frac{\mu}{2}\left(\bm{C}:\bm{I}-3\right)-\mu\log\left(\sqrt{\det\bm{C}}\right)+\frac{\chi}{2}\log^{2}\left(\sqrt{\det\bm{C}}\right)$
with $\mu=\frac{E}{2(1+\nu)}$ and $\chi=\frac{E\nu}{(1+\nu)(1-2\nu)}$
being constitutive properties.

\subsection*{Patch test satisfaction}

We employ a corrected penalty so that the contact patch test is satisfied in most cases. This is an important subject for convergence of computational contact solutions and has been addressed here with a similar solution to the one discussed by Zavarise and co-workers \cite{zavarise2009,zavarise2009b}.

\begin{figure}
\centering
\includegraphics[clip,width=0.8\textwidth]{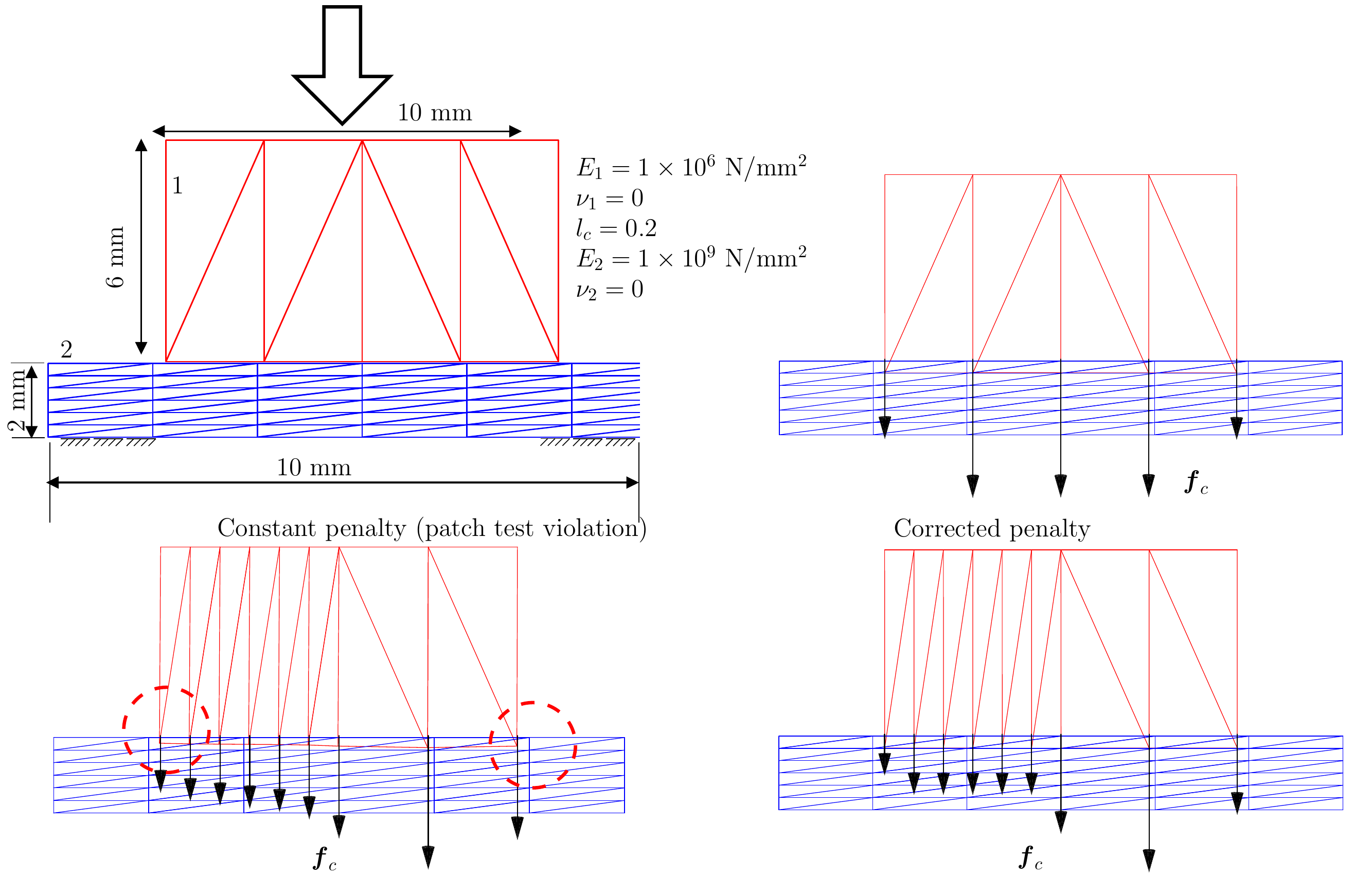}
\caption{Patch test satisfaction by using a corrected penalty.}
\label{fig:pts}
\end{figure}

We remark that this is not a general solution and in some cases, our formulation may fail to pass the patch test. Figure \ref{fig:pts} shows the effect of using a penalty weighted by the edge projection, see \cite{zavarise2009b}. However, this is not an universal solution. 

\subsection{Two-dimensional compression}
\begin{figure}
\centering
\includegraphics[clip,width=0.8\textwidth]{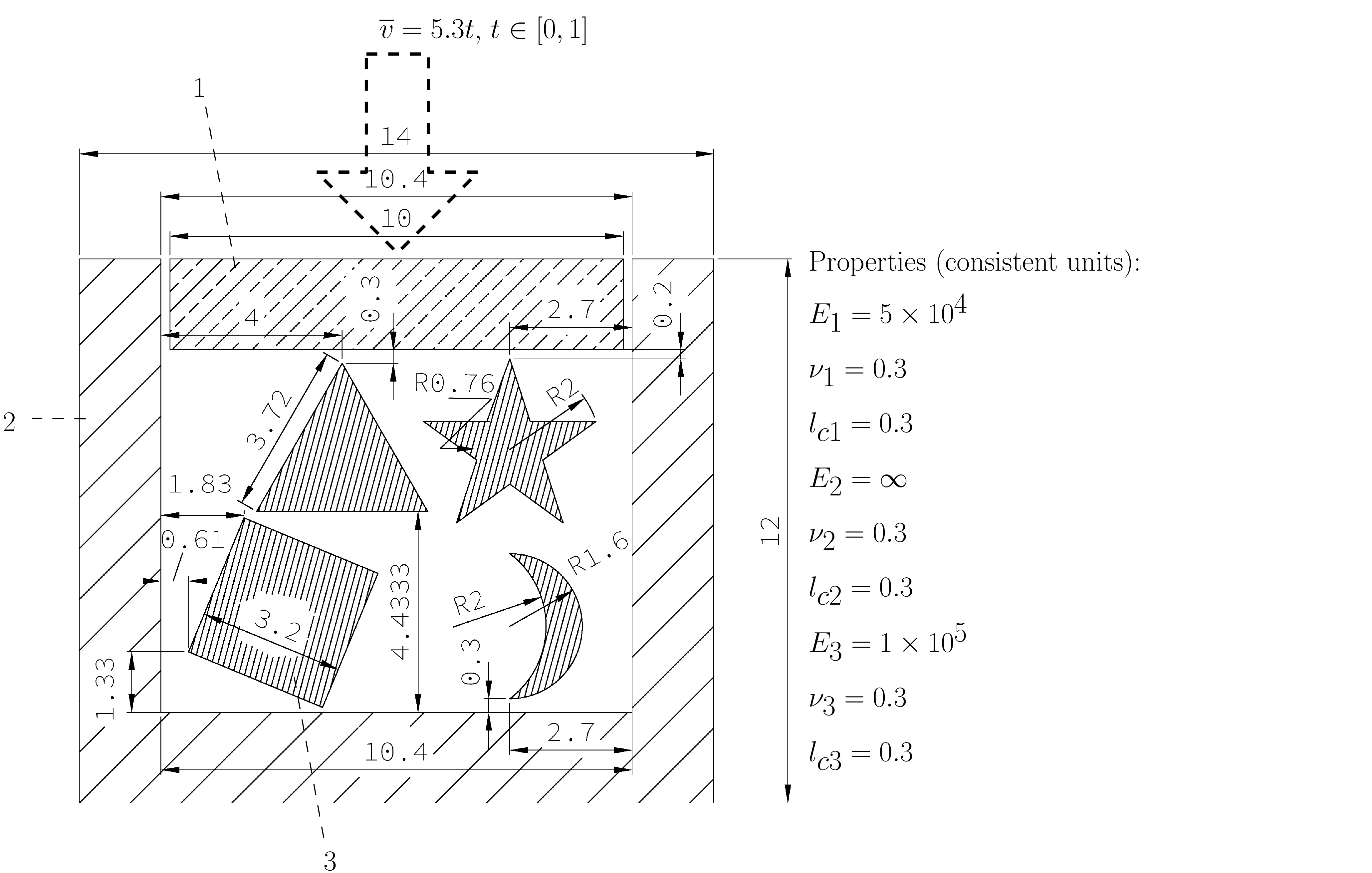}
\caption{Two-dimensional verification problem. Consistent units are used.}
\label{fig:Illustrative-problem-(2D)}
\end{figure}

We begin with a quasi-static two-dimensional test, as shown in Figure~\ref{fig:Illustrative-problem-(2D)},
using three-node triangles. This test consists of a compression of
four polygons (identified as part 3) in Figure \ref{fig:Illustrative-problem-(2D)}
by a deformable rectangular punch (part 1 in the same Figure). The
`U' (part 2) is considered rigid but still participates in the
approximate distance function (ADF) calculation. To avoid an underdetermined
system, a small damping term is used, specifically $40$ units with
$\mathsf{L}^{-2}\mathsf{M}\mathsf{T}^{-1}$ ISQ dimensions. Algorithm
\ref{alg:Contact-algorithm-based} is adopted with a pseudo-time increment
of $\Delta t=0.003$ for $t\in[0,1]$.

\begin{figure}
\begin{centering}
\begin{tabular}{cccc}
\multicolumn{4}{c}{$\overline{v}=0.00$}\tabularnewline
\multicolumn{4}{c}{\includegraphics[width=0.3609\textwidth]{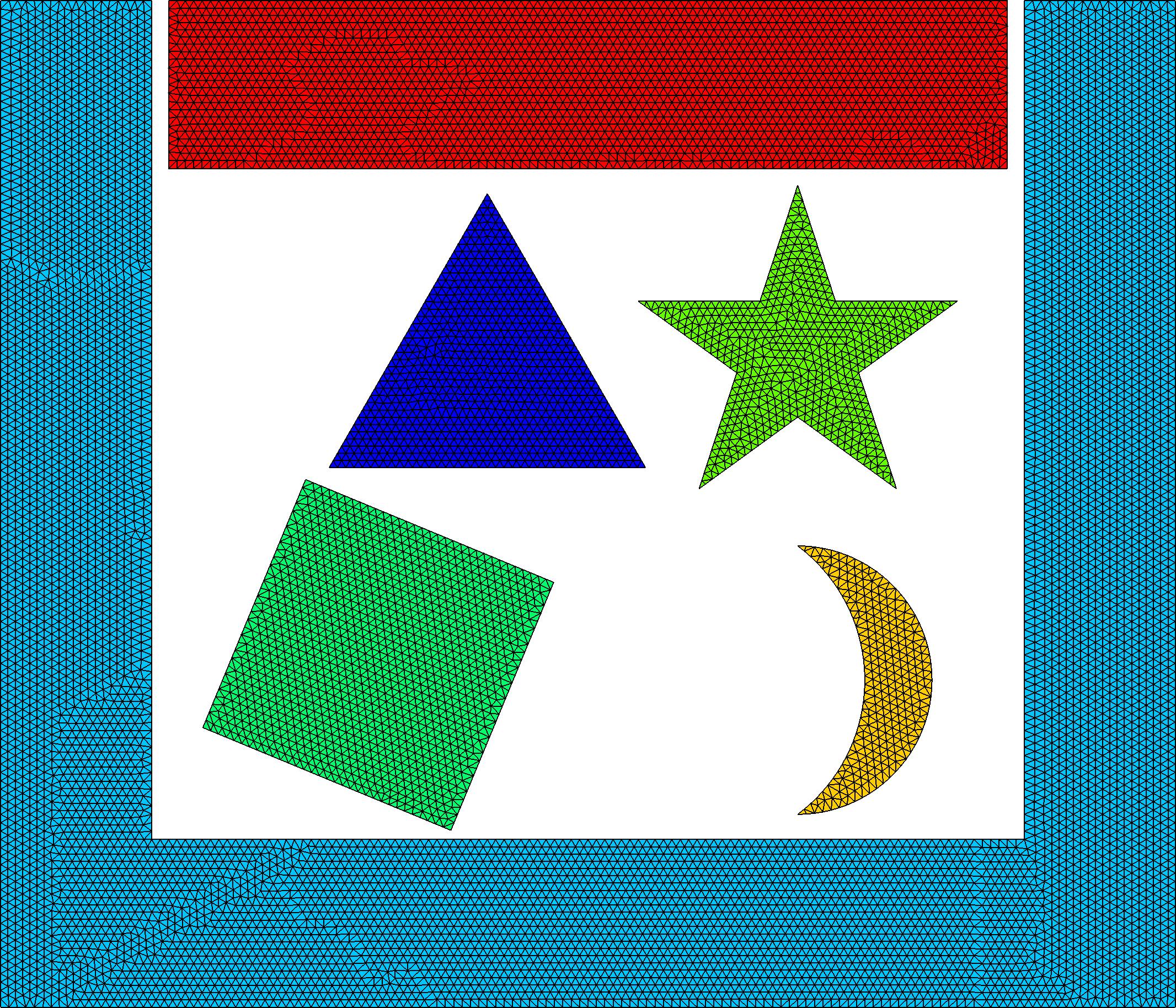}\hspace{1cm}\includegraphics[width=0.4941\textwidth]{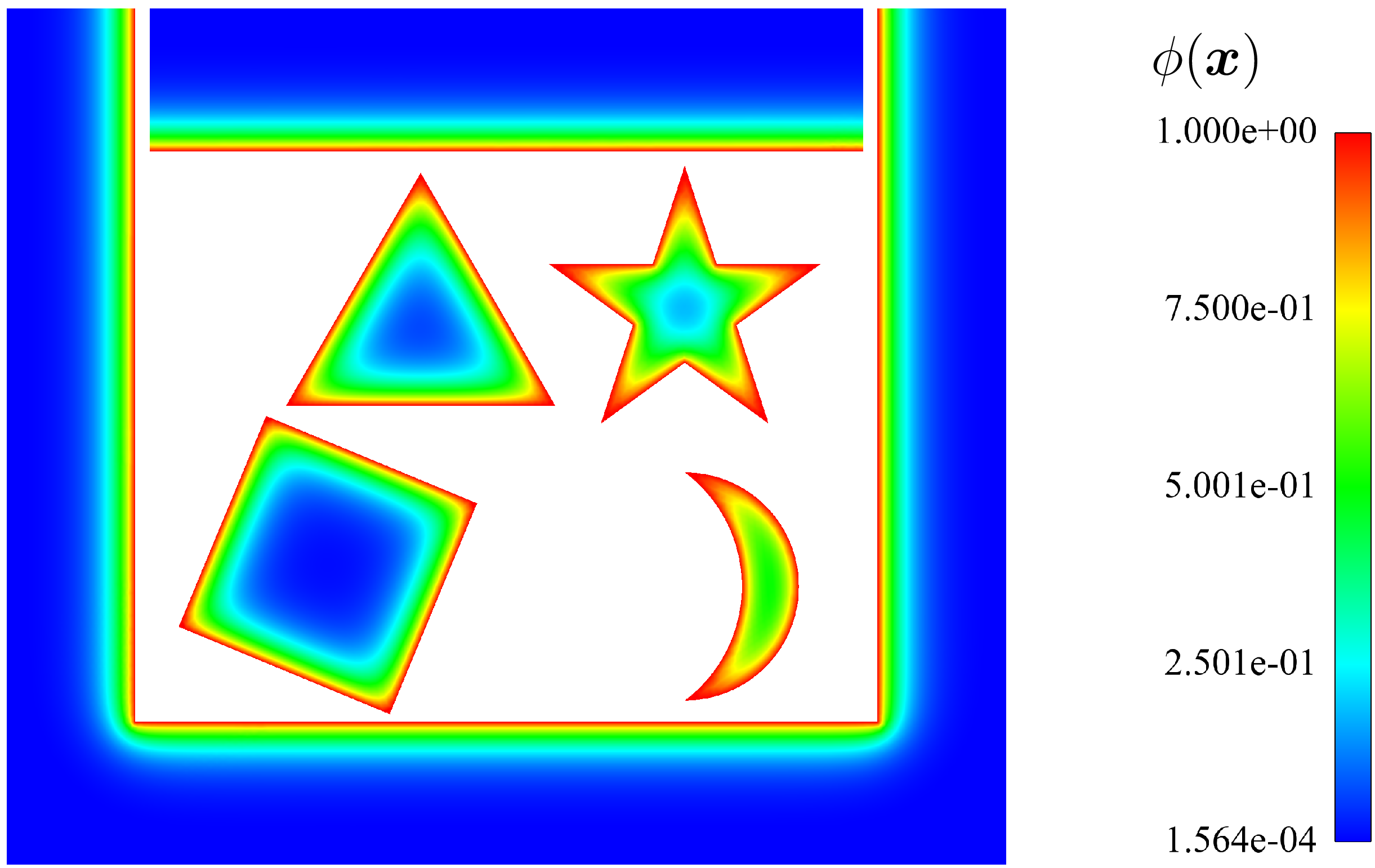}}\tabularnewline
 & $\overline{v}=3.24$ & $\overline{v}=5.30$ & $\overline{v}=5.30$, $\phi$\tabularnewline
{\Large{}\rotatebox{90}{$h=0.020$}} & \includegraphics[width=0.27702\textwidth]{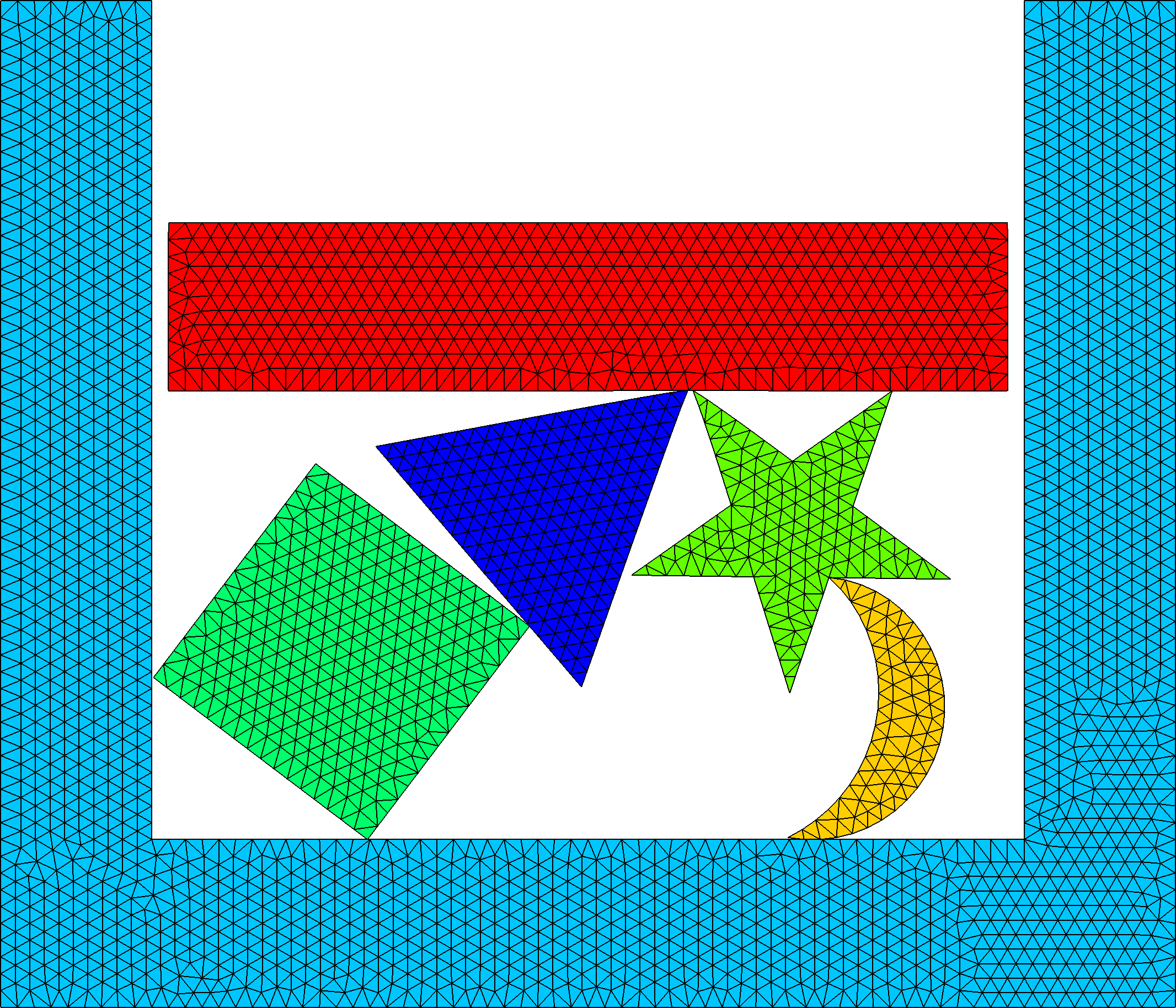} & \includegraphics[width=0.27702\textwidth]{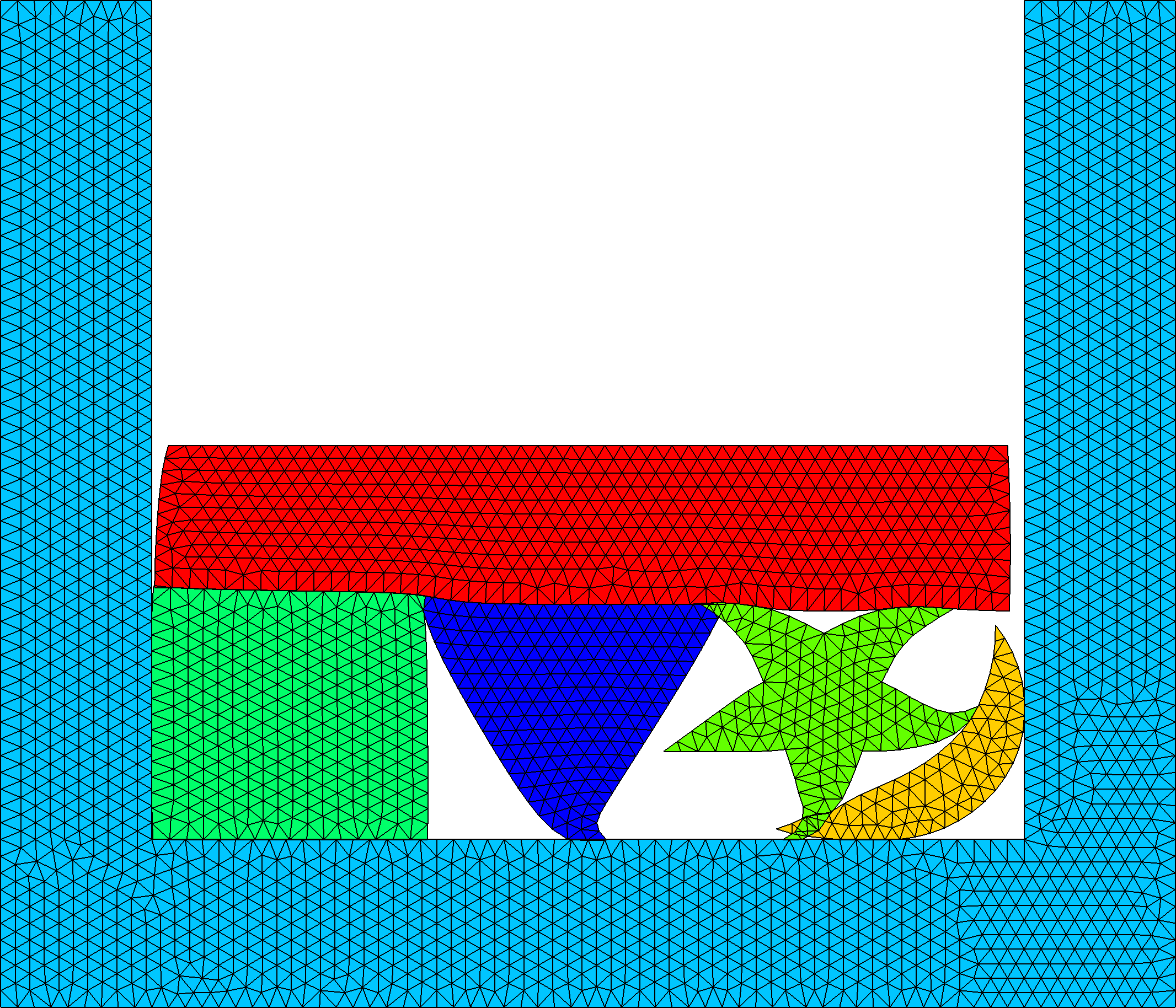} & \includegraphics[width=0.3807\textwidth]{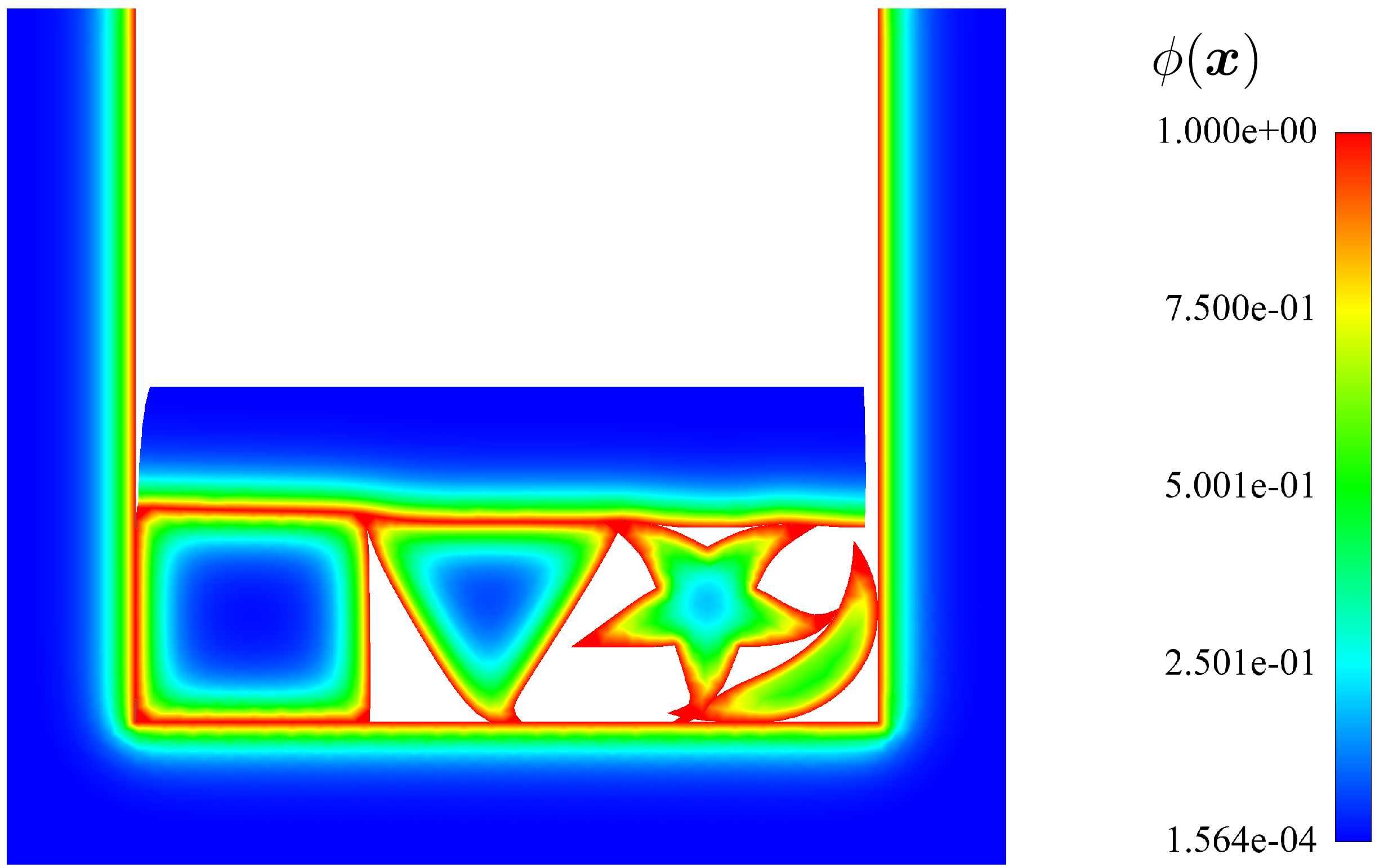}\tabularnewline
{\Large{}\rotatebox{90}{$h=0.015$}} & \includegraphics[width=0.27702\textwidth]{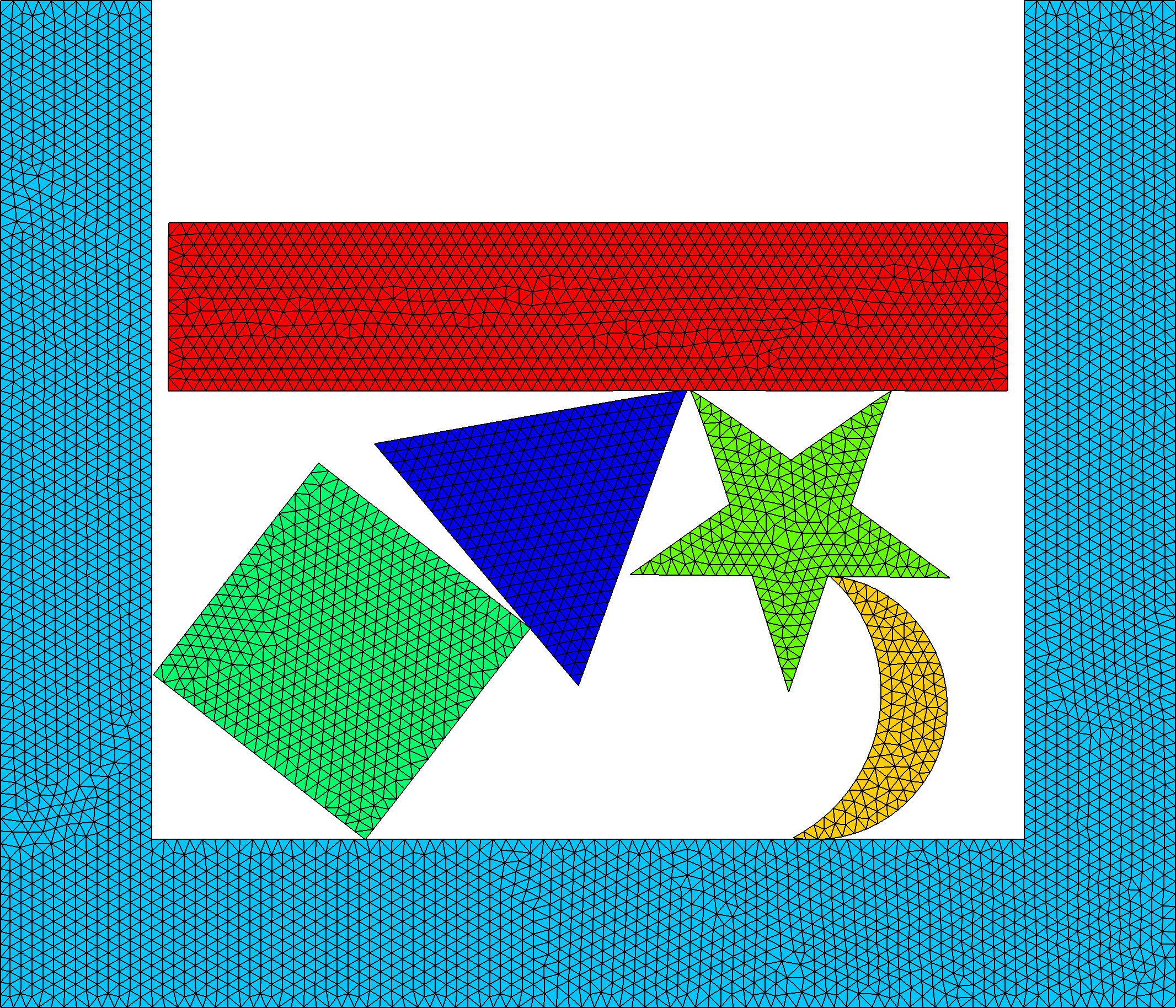} & \includegraphics[width=0.27702\textwidth]{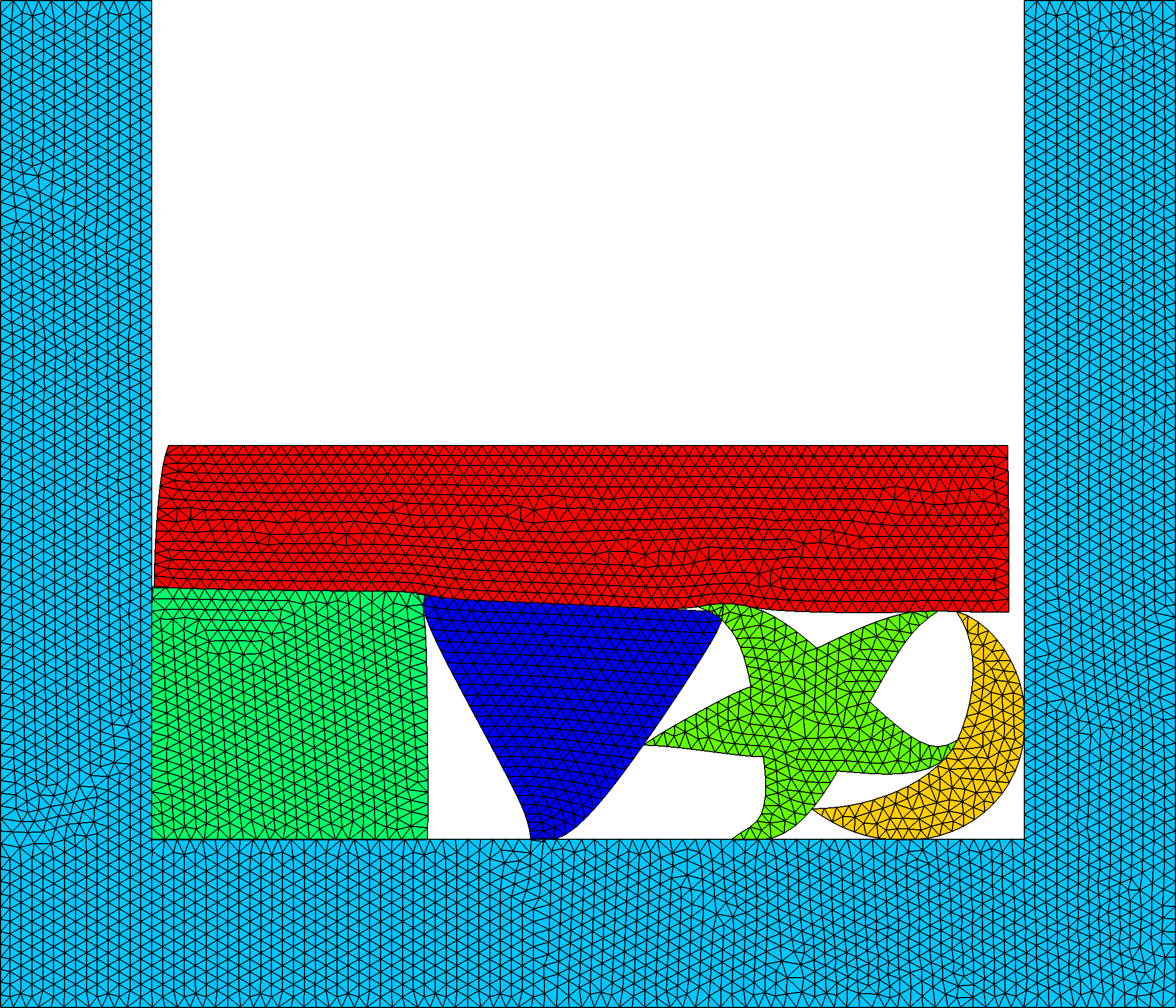} & \includegraphics[width=0.3807\textwidth]{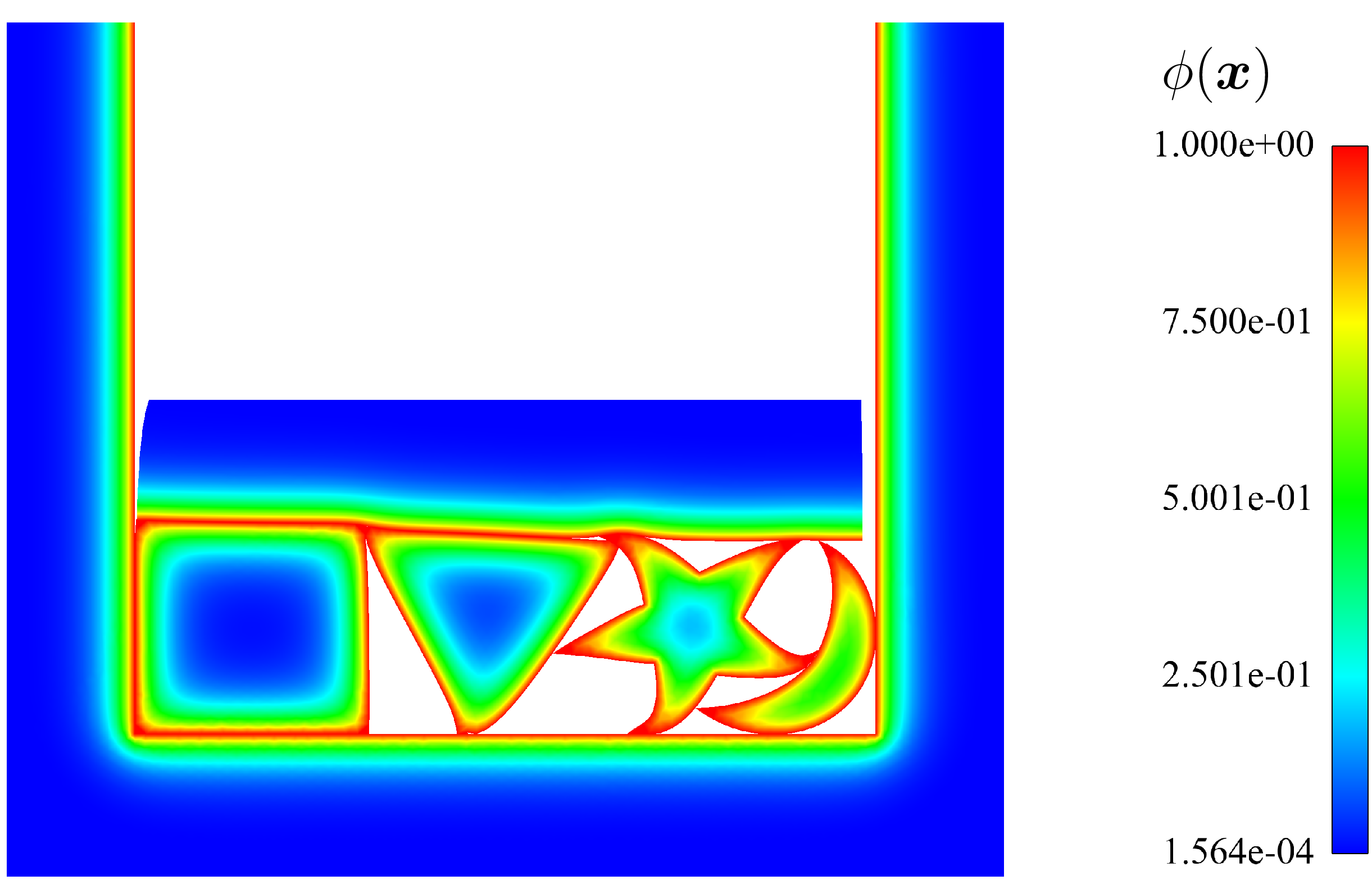}\tabularnewline
{\Large{}\rotatebox{90}{$h=0.010$}} & \includegraphics[width=0.27702\textwidth]{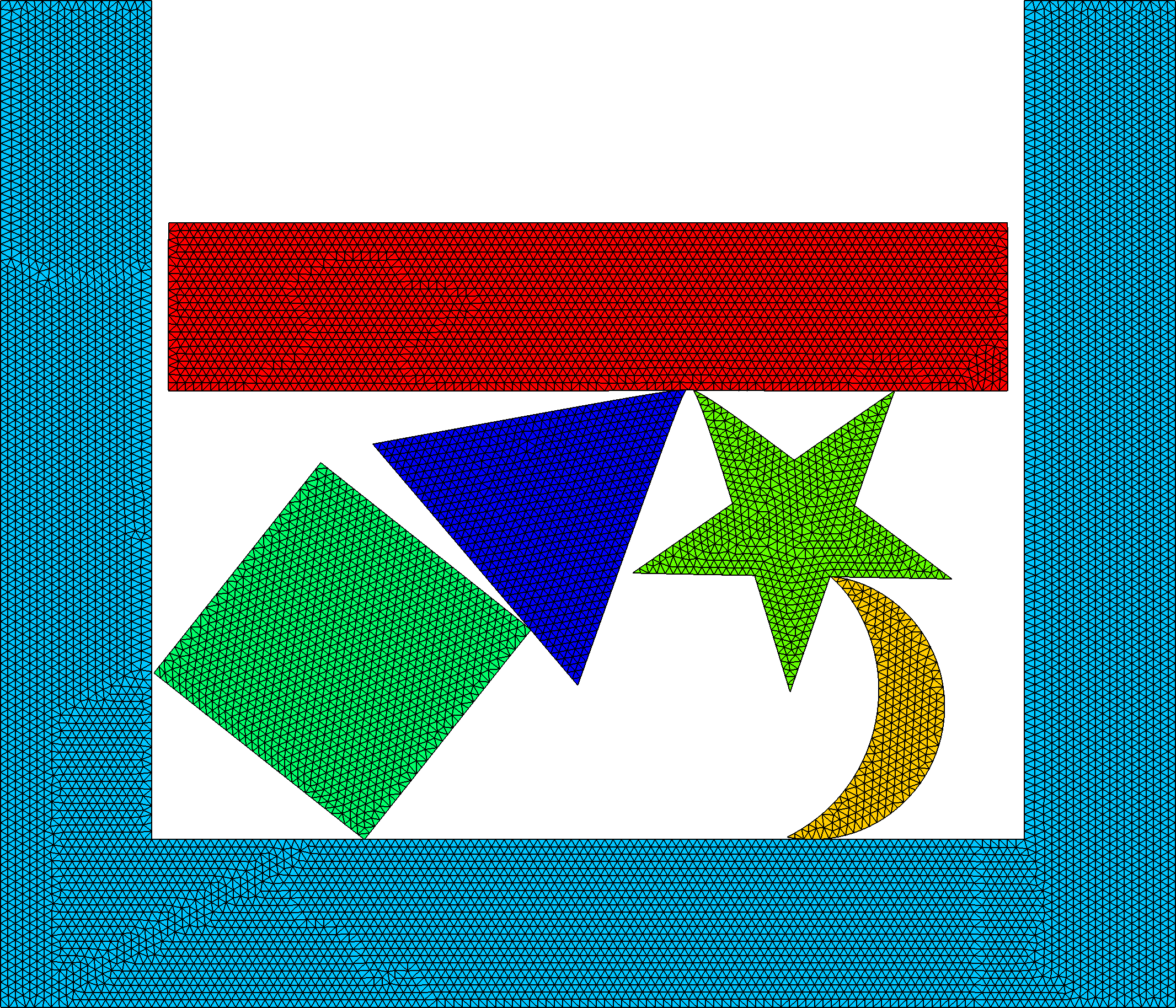} & \includegraphics[width=0.27702\textwidth]{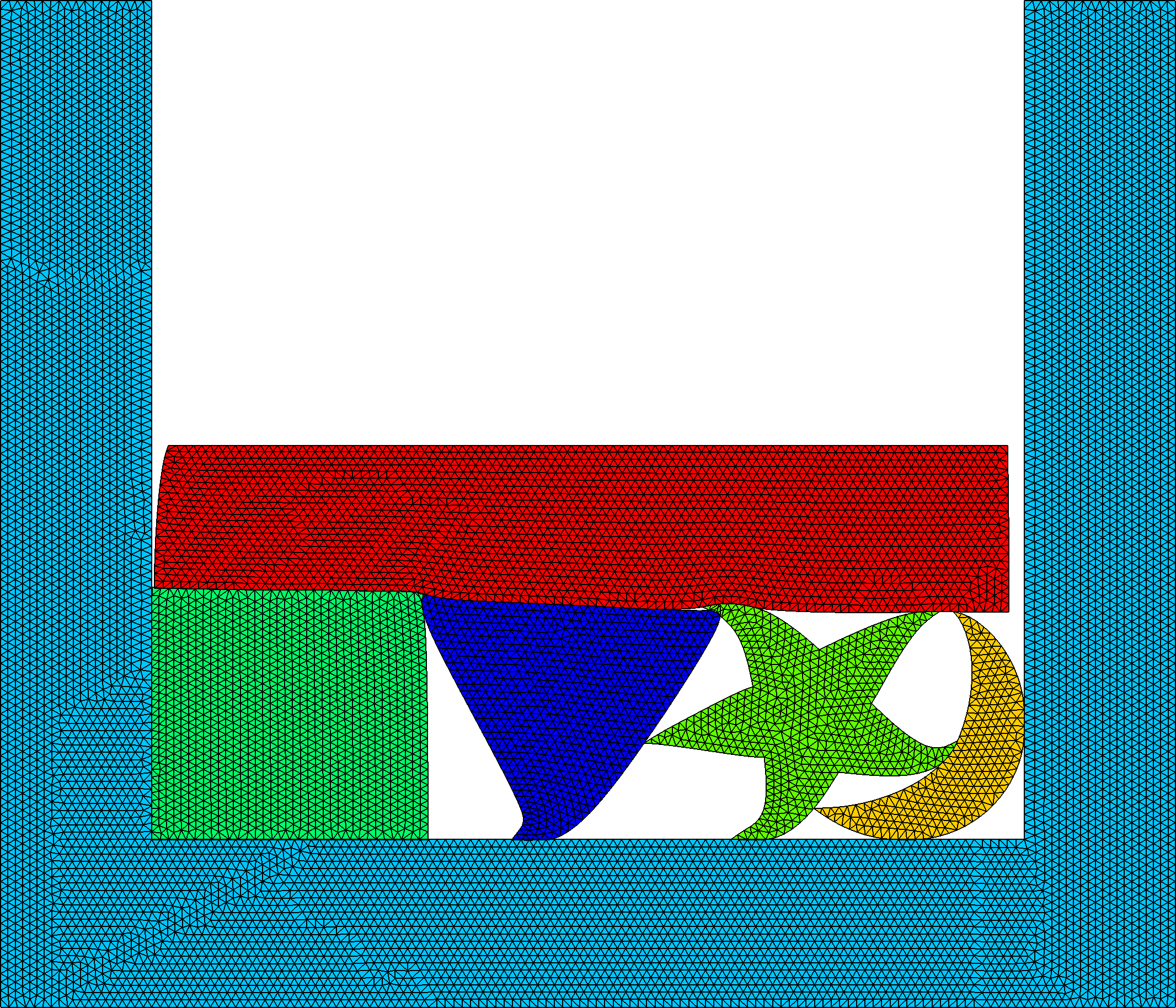} & \includegraphics[width=0.3807\textwidth]{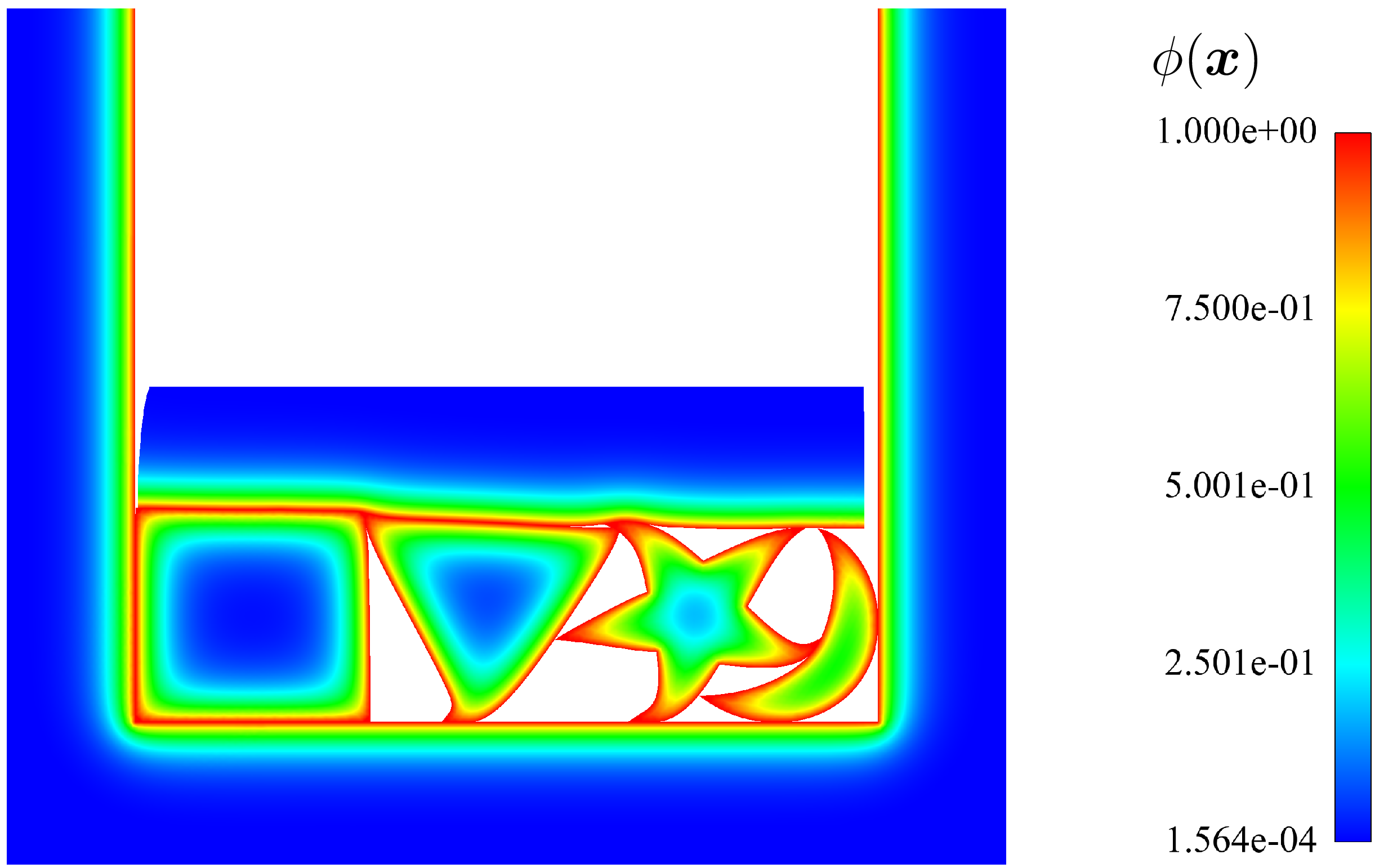}\tabularnewline
\end{tabular}
\par\end{centering}
\caption{\label{fig:Distance-interference}Two-dimensional compression: sequence of deformed meshes and contour plot of $\phi(\bm{x})$.}
\end{figure}

\begin{figure}
\centering
\includegraphics[width=0.9\textwidth]{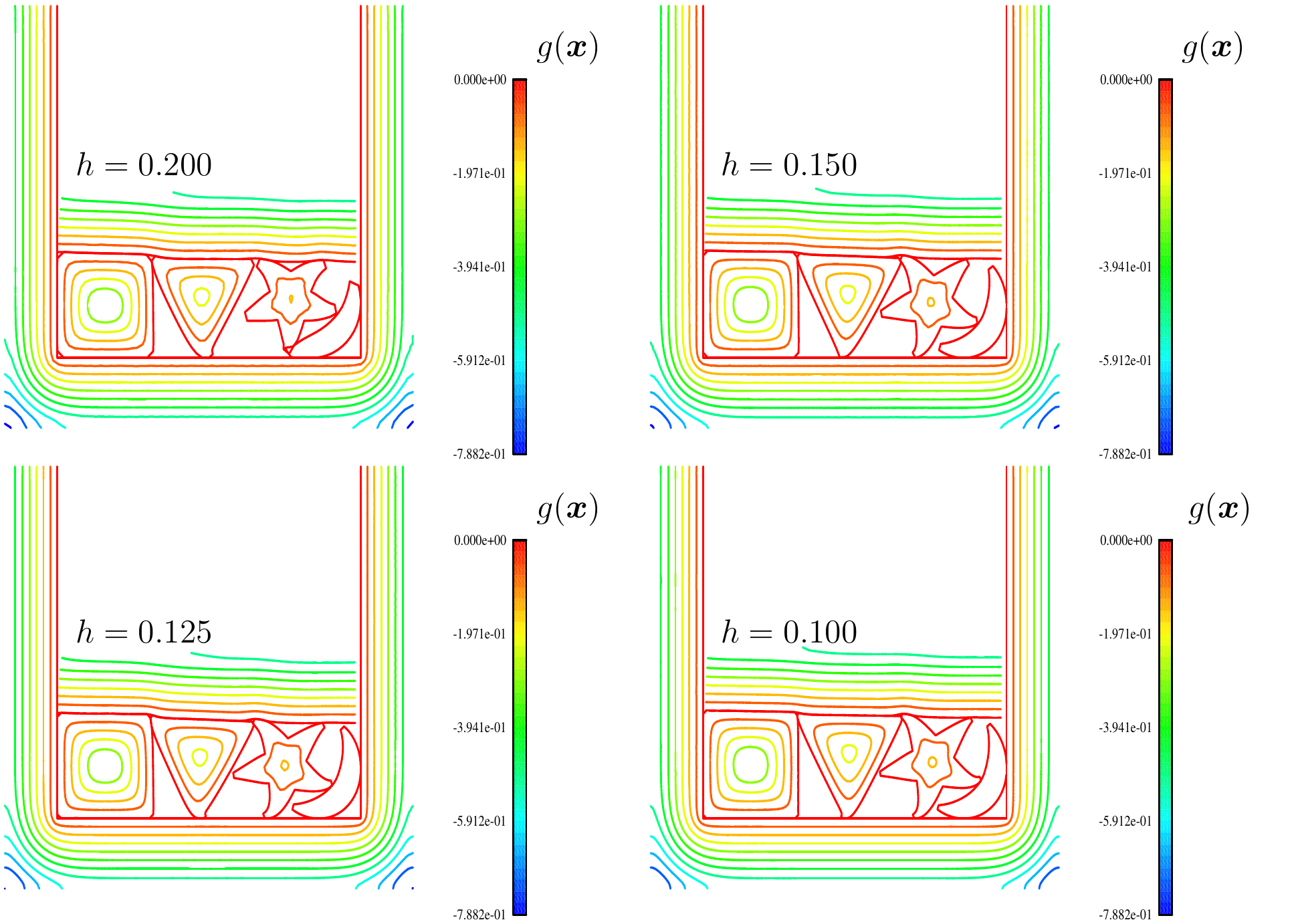}
\caption{\label{fig:col}Gap ($g(\bm{x})$) contour lines for the four meshes using $l_{c}=0.3$ consistent units.}
\end{figure}

For $h=0.020$, $h=0.015$ and $h=0.010$, Figure \ref{fig:Distance-interference} shows a sequence of deformed meshes and the contour plot of $\phi(\bm{x})$. The robustness of the algorithm is excellent, at the cost of some interference for coarse meshes. To further inspect the interference, the contour lines for $g(\bm{x})$ are shown in Figure \ref{fig:col}. We note that coarser meshes produce smoothed-out vertex representation, which causes the interference displayed in Figure \ref{fig:Distance-interference}. Note that $g(\bm{x})$ is determined from 
$\phi(\bm{x})$.
%
%

\begin{figure}
\begin{centering}
\subfloat[$h=0.020$]{\begin{centering}
\includegraphics[width=0.45\textwidth]{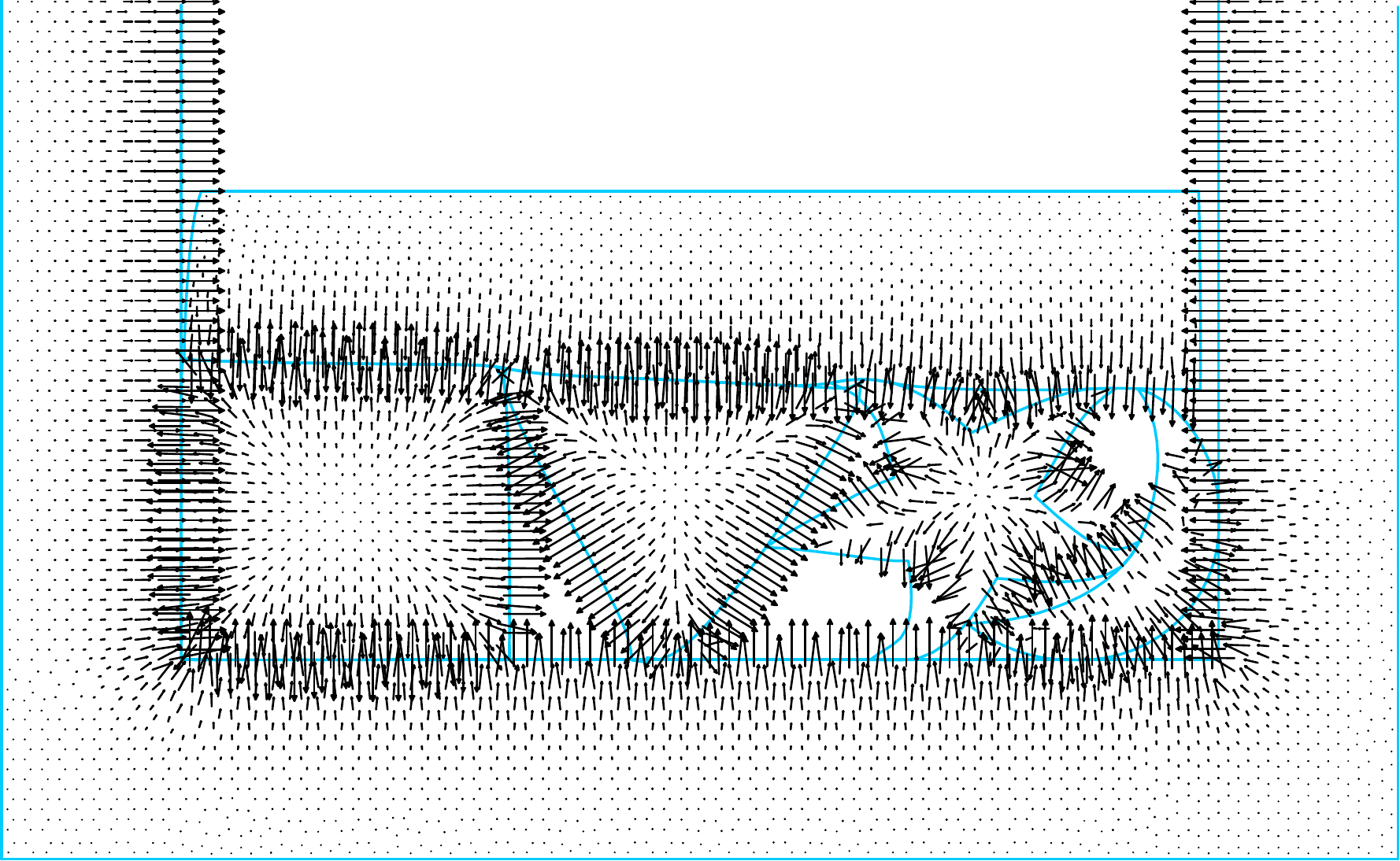}
\par\end{centering}
}
\par\end{centering}
\begin{centering}
\subfloat[$h=0.015$]{\begin{centering}
\includegraphics[width=0.45\textwidth]{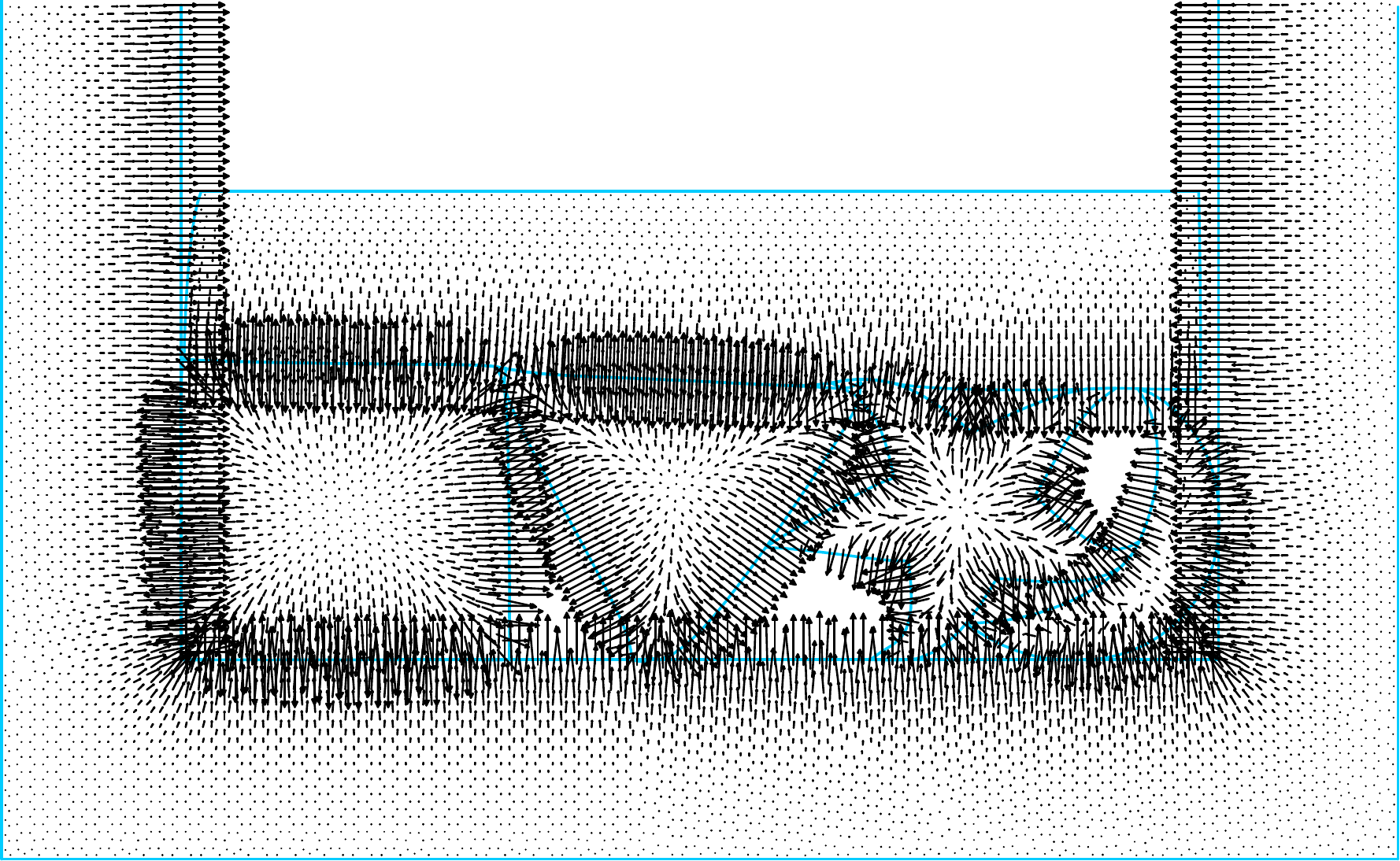}
\par\end{centering}
}
\par\end{centering}
\begin{centering}
\subfloat[$h=0.010$]{\begin{centering}
\includegraphics[width=0.45\textwidth]{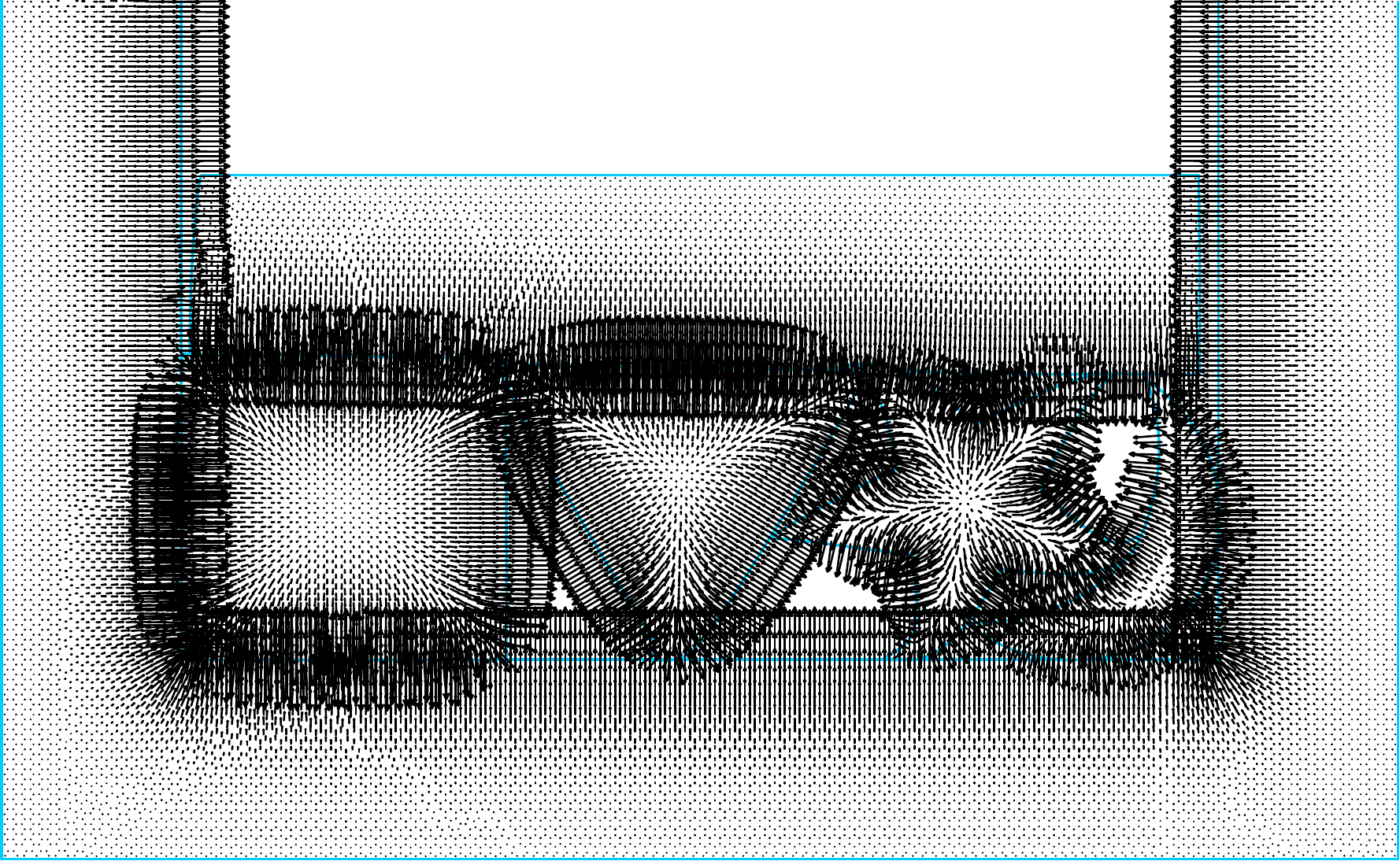}
\par\end{centering}
}
\par\end{centering}
\caption{\label{fig:Directions-obtained-as}Directions obtained as $\nabla\phi\left(\bm{\xi}\right)$
for $h=0.020,$ $0.015$ and $0.010$.}
\end{figure}

Using the gradient of $\phi\left(\bm{x}\right)$, the contact direction
is obtained for $h=0.02$ as shown in Figure~\ref{fig:Directions-obtained-as}.
We can observe that for the star-shaped figure, vertices are poorly
represented, since small gradients are present due to uniformity of
$\phi\left(\bm{x}\right)$ in these regions. The effect of mesh refinement
is clear, with finer meshes producing a sharper growth of reaction
when all four objects are in contact with each other. In contrast,
the effect of the characteristic length $l_{c}$ is not noticeable.

\begin{figure}
\begin{centering}
\includegraphics[width=0.8\textwidth]{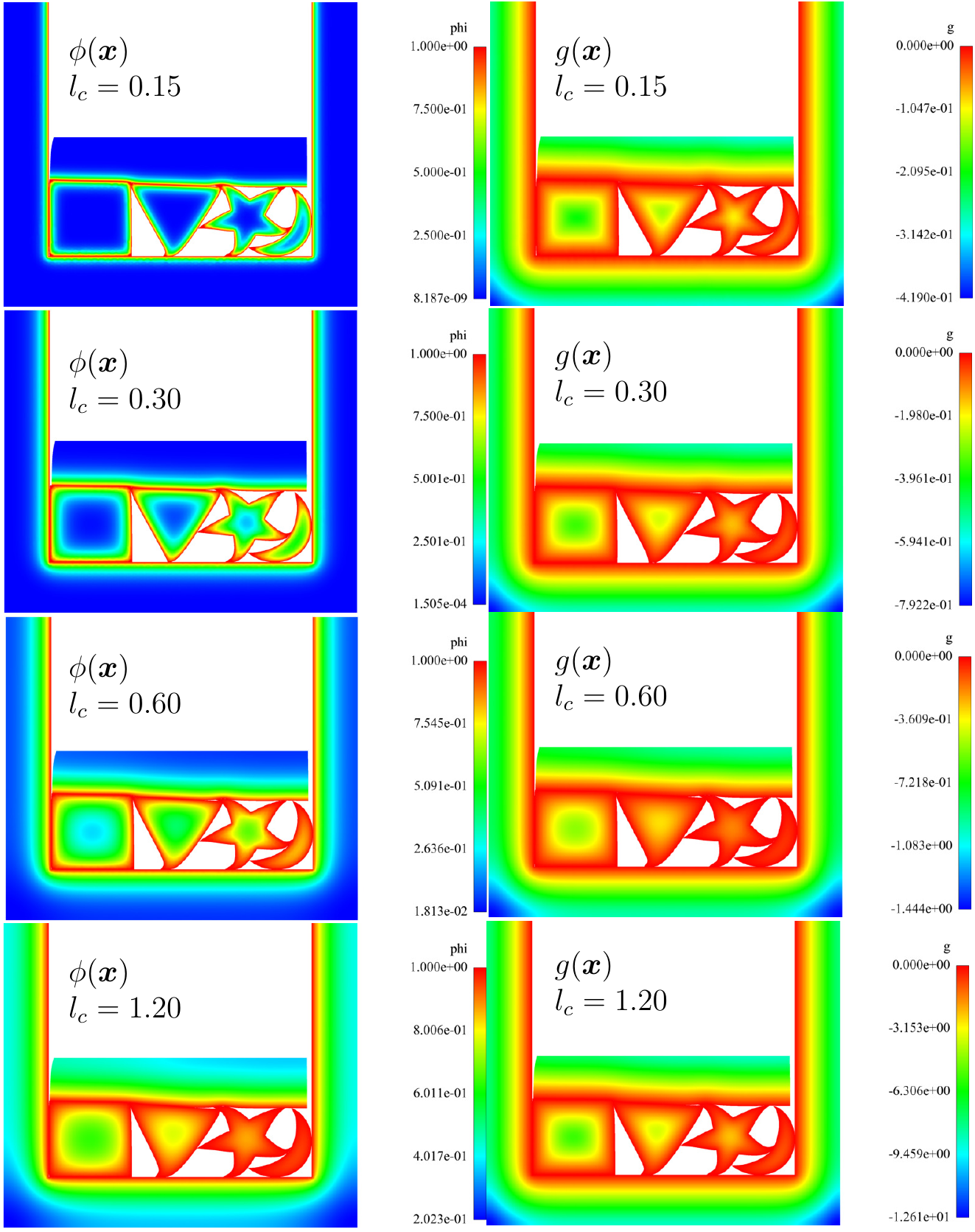}
\par\end{centering}
\caption{\label{fig:lcphig}Effect of $l_{c}$ in the form of $\phi(\bm{x})$ and $g(\bm{x})$.}
\end{figure}

\begin{figure}
\begin{centering}
\subfloat[]{\begin{centering}
\includegraphics[width=0.50\textwidth]{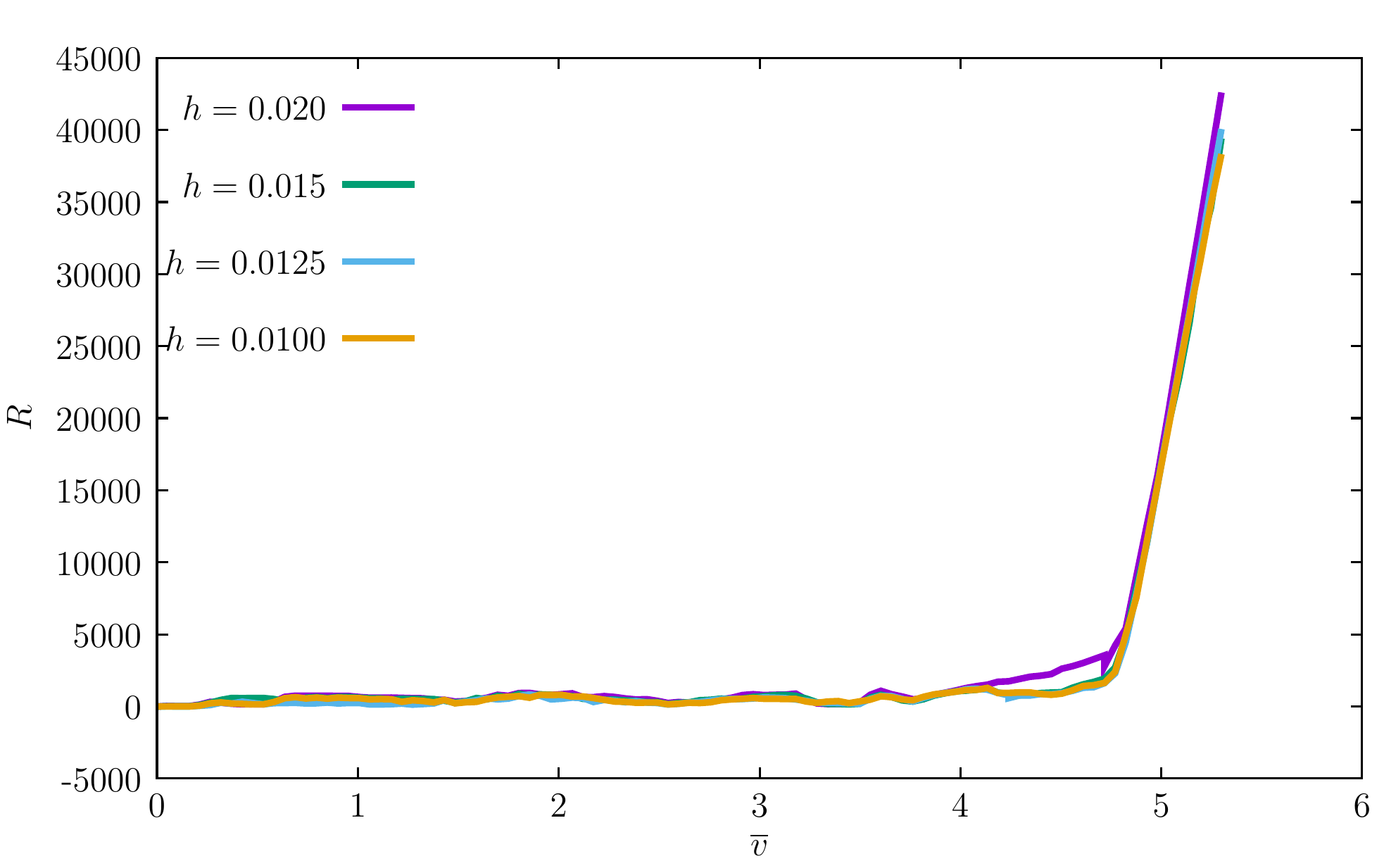}
\par\end{centering}
}\subfloat[]{\begin{centering}
\includegraphics[width=0.50\textwidth]{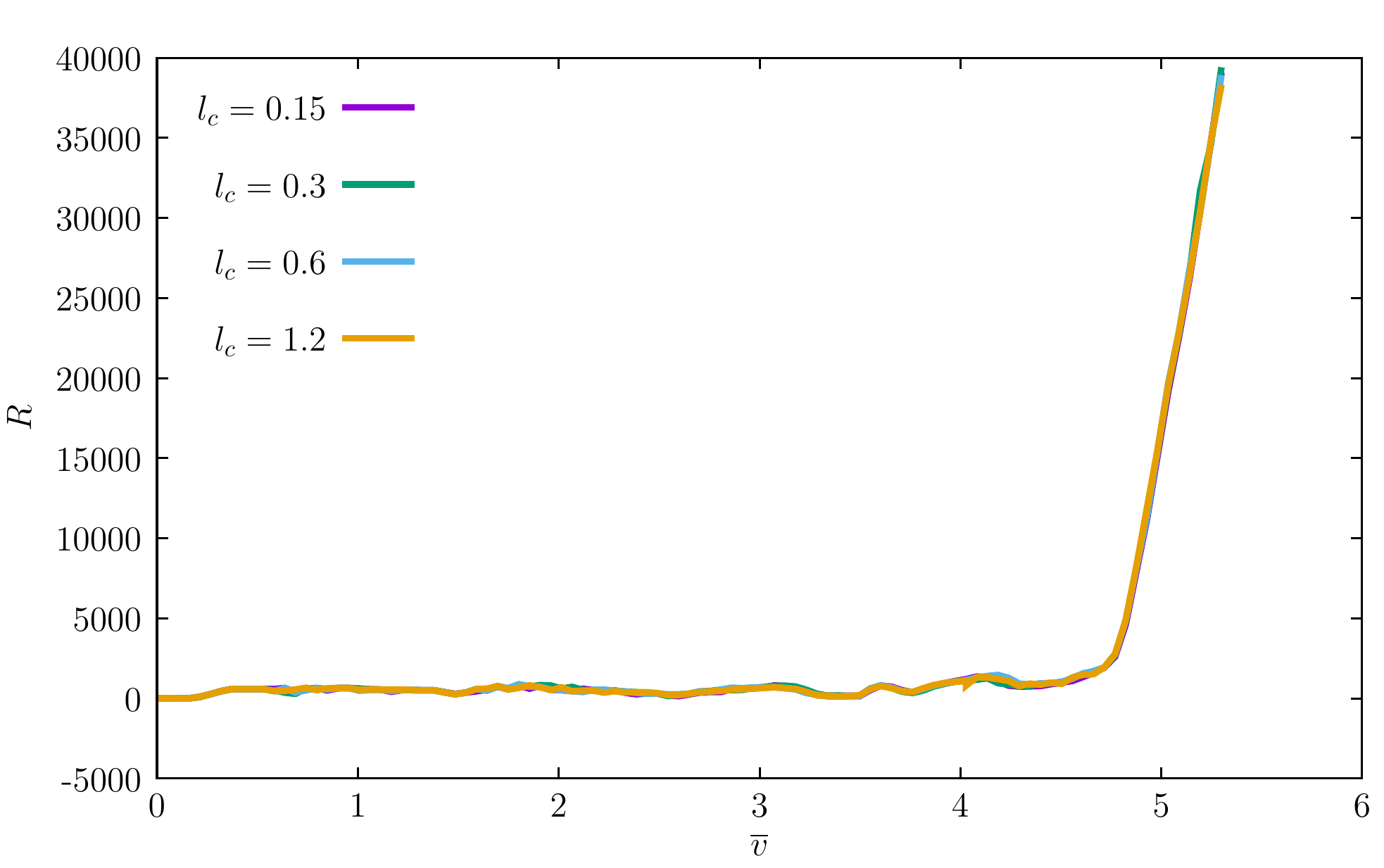}
\par\end{centering}
}
\par\end{centering}
\caption{\label{fig:Effect-of-}Effect of $h$ in (a) and effect of $l_{c}$ in
(b) in the displacement/reaction behavior.}
\end{figure}

In terms of the effect of $l_c$ on the fields $\phi(\bm{x})$ and $g(\bm{x})$, Figure~\ref{fig:lcphig} shows that information, although $\phi(\bm{x})$ is strongly affected by the length parameter, $g(\bm{x})$ shows very similar spatial distributions although different peaks. Effects of $h$ and $l_c$ in the displacement/reaction behavior is shown in Figure \ref{fig:Effect-of-}. The mesh size $h$ has a marked effect on the results up to $h=0.0125$, the effect of $l_c$ is much weaker.

\FloatBarrier

\subsection{Three-dimensional compression}
\begin{figure}
\begin{centering}
\subfloat[]{\begin{centering}
\includegraphics[width=0.7\textwidth]{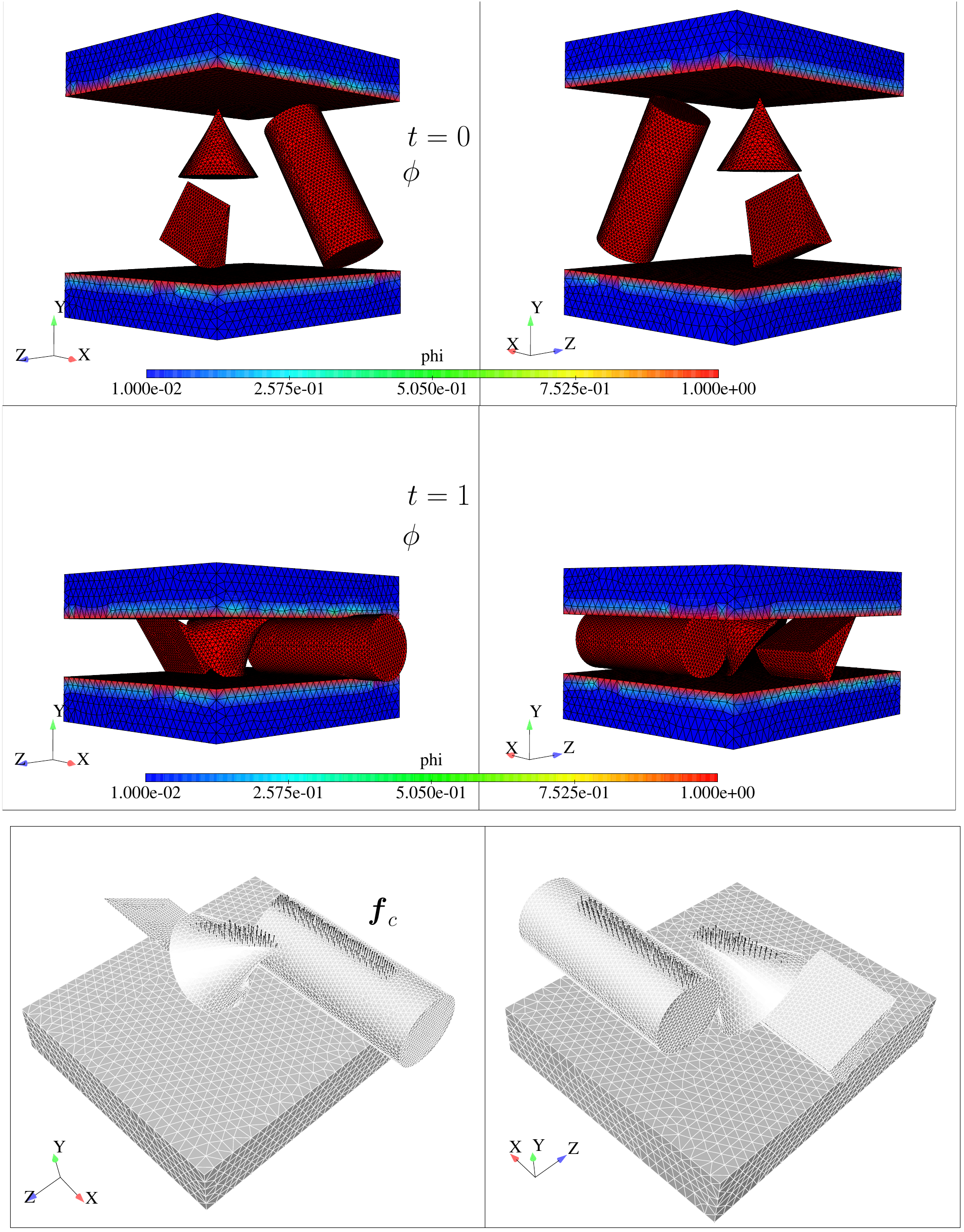}
\par\end{centering}
}
\par\end{centering}
\begin{centering}

\subfloat[]{\begin{centering}
\includegraphics[width=0.6\textwidth]{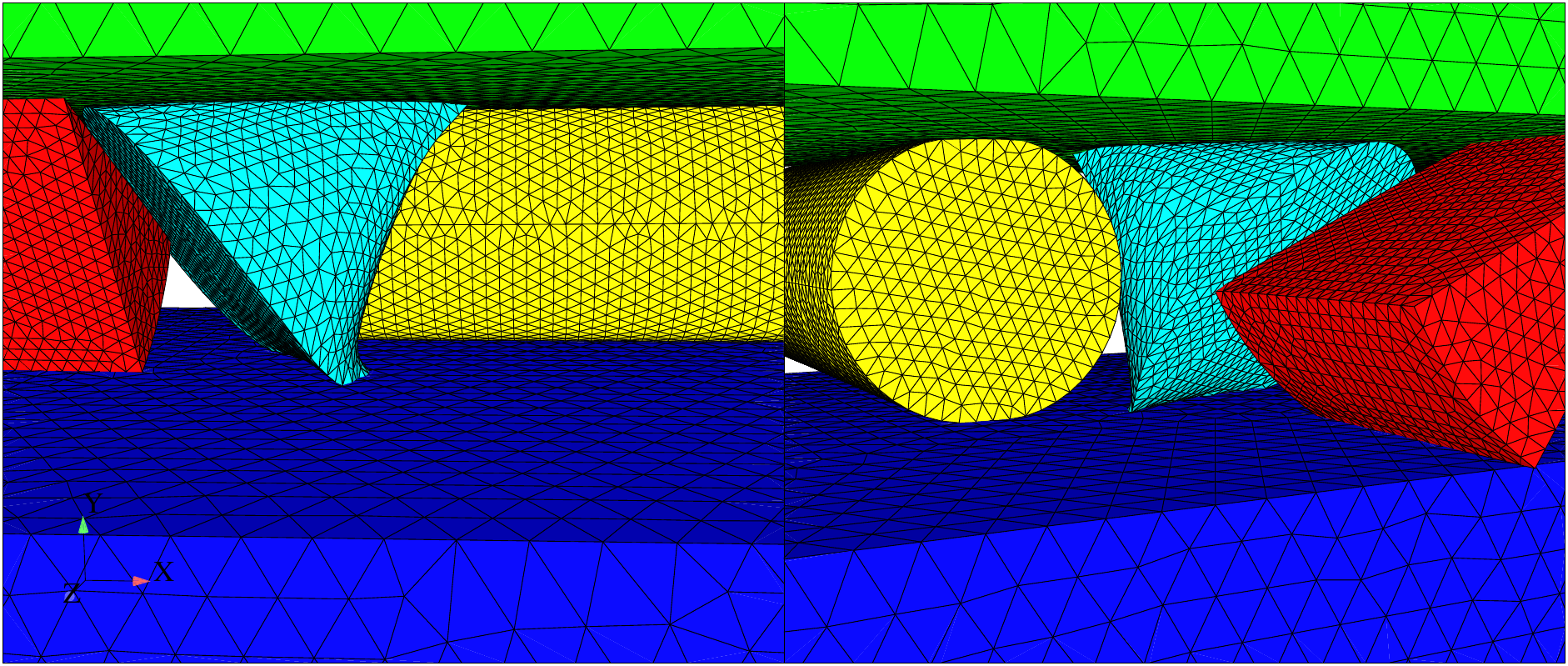}
\par\end{centering}
}
\par\end{centering}
\caption{Data for the three-dimensional compression
test. (a) Undeformed and deformed configurations ($h=0.025$). For geometric
files, see~\cite{Areias2022cgithub}.
(b) Detail of contact of cone with lower block and with the wedge. Each object is identified by a different color.}
\label{fig:Three-dimensional-...}
\end{figure}

In three dimensions, the algorithm is in essence the same. Compared
with geometric-based algorithms, it significantly reduces the coding
for the treatment of particular cases (node-to-vertex, node-to-edge
in both convex and non-convex arrangements). The determination of
coordinates for each incident node is now performed on each tetrahedra,
but the remaining tasks remain unaltered. We test the geometry shown
in Figure \ref{fig:Three-dimensional-...} with the following objectives:
\begin{itemize}
\item Assess the extension to the 3D case. Geometrical nonsmoothness
is introduced with a cone and a wedge.
\item Quantify interference as a function of $l_{c}$ and $\kappa$ as well
as the strain energy evolution.
\end{itemize}
Deformed configurations and contour plots of $\phi\left(\bm{x}\right)$
for this problem are presented in Figure \ref{fig:Three-dimensional-...},
and the corresponding CAD files are available on \texttt{Github}~\cite{Areias2022cgithub}.
A cylinder, a cone and a wedge are compressed between two blocks.
Dimensions of the upper and lower blocks are $10\times12\times2$
consistent units (the upper block is deformable whereas the lower
block is rigid) and initial distance between blocks is $8$ consistent
units. Length and diameter of the cylinder are $7.15$
and $2.86$ (consistent units), respectively. The cone has a height of $3.27$ consistent
units and a radius of $1.87$. Finally, the wedge has a width of $3.2$,
a radius of $3.2$ and a swept angle of $30$ degrees. A compressible
Neo-Hookean law is adopted with the following properties:
\begin{itemize}
\item Blocks: $E=5\times10^{4}$ and $\nu=0.3$.
\item Cylinder, cone and wedge:~$E=1\times10^{5}$ and $\nu=0.3$.
\end{itemize}

\begin{figure}
\centering
\subfloat[]{\begin{centering}
\includegraphics[width=0.50\textwidth]{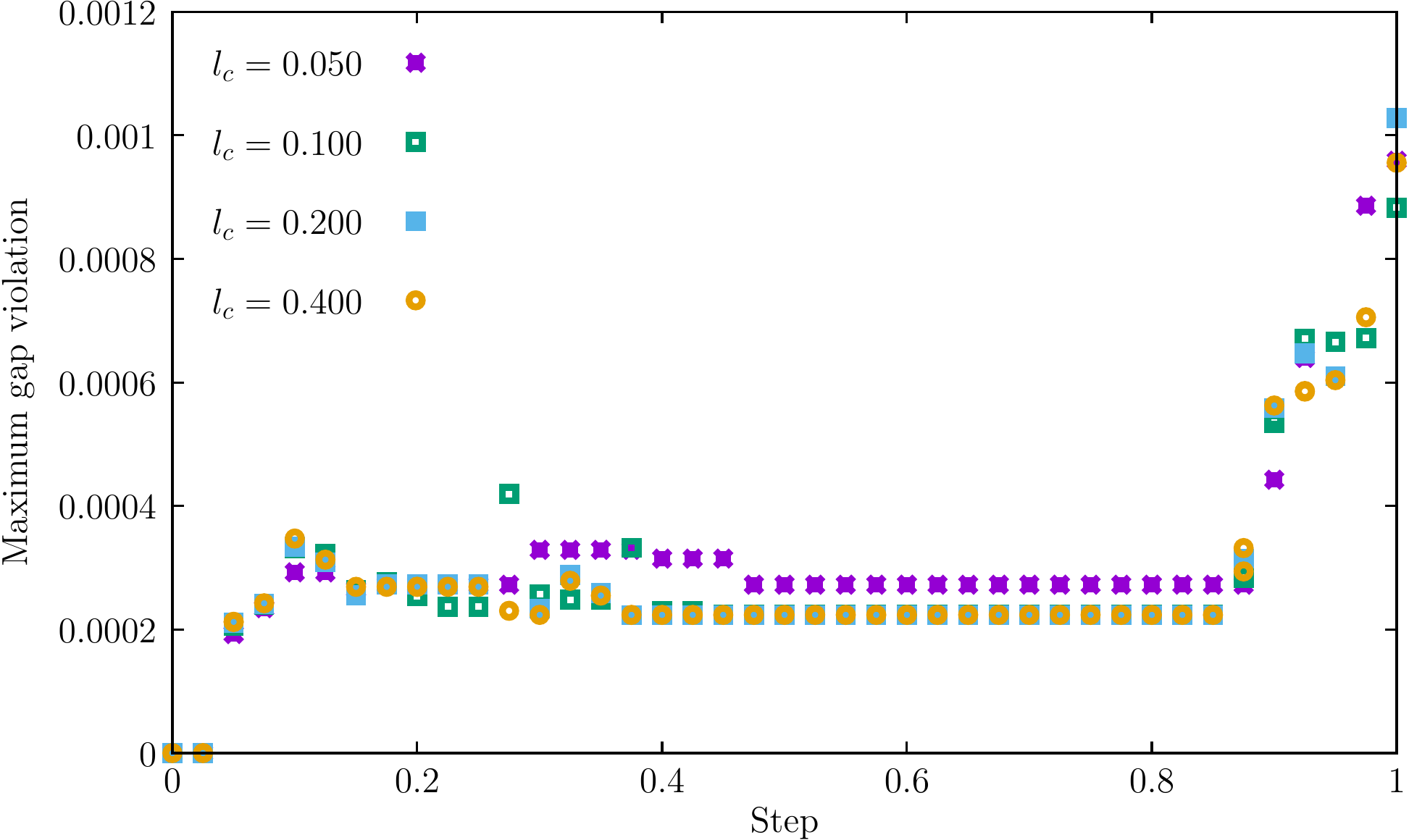}
\par\end{centering}
}\subfloat[]{\begin{centering}
\includegraphics[width=0.50\textwidth]{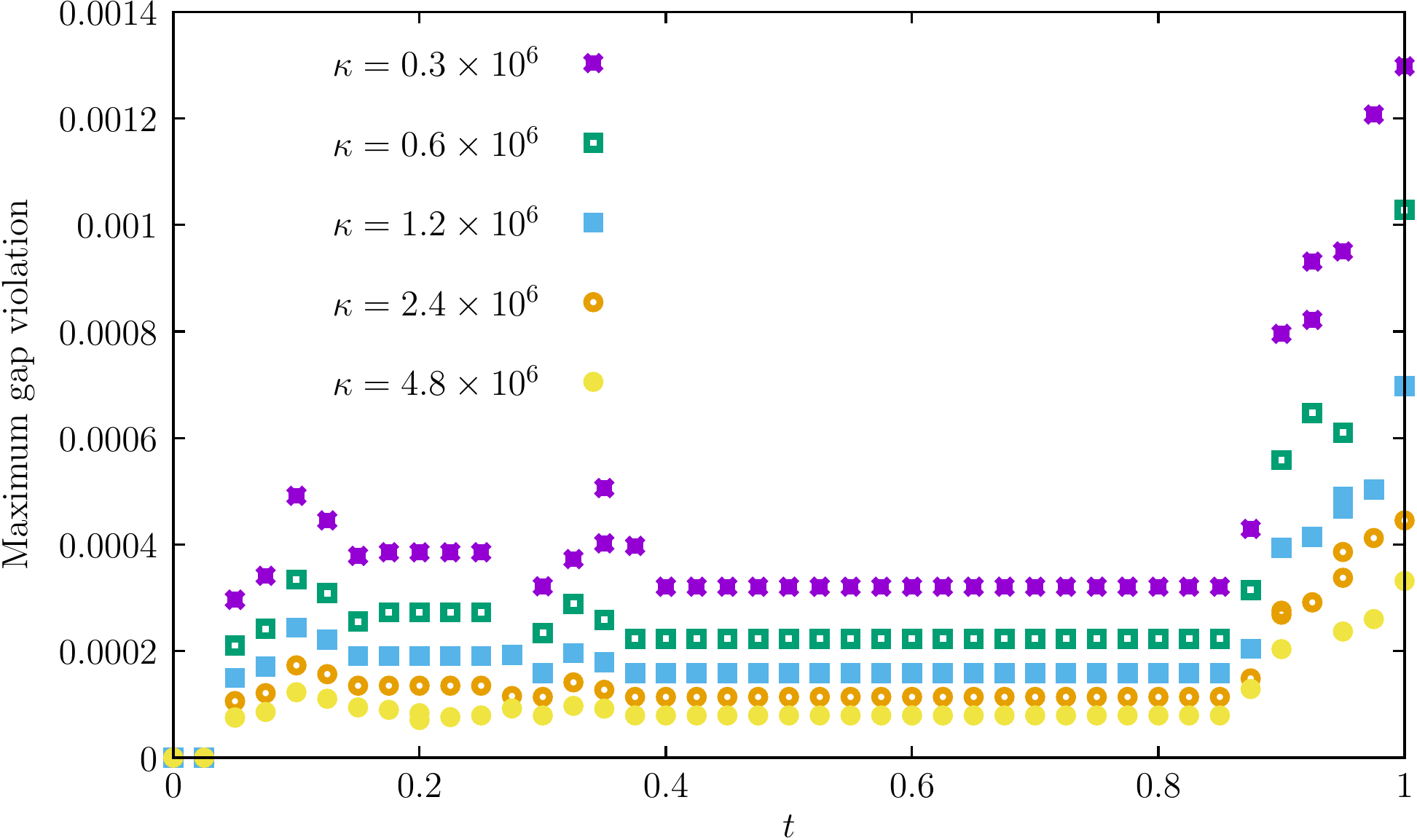}
\par\end{centering}
}
\caption{(a) Effect of $l_{c}$ on the maximum gap evolution
over pseudotime $(\kappa=0.6\times10^{6}$, $h=0.030)$ and
(b) Effect of $\kappa$ on the maximum gap evolution over pseudotime 
$(l_{c}=0.2$).}
\label{fig:Effect-of--1}
\end{figure}

\begin{figure}
\centering
\subfloat[]{\begin{centering}
\includegraphics[width=0.50\textwidth]{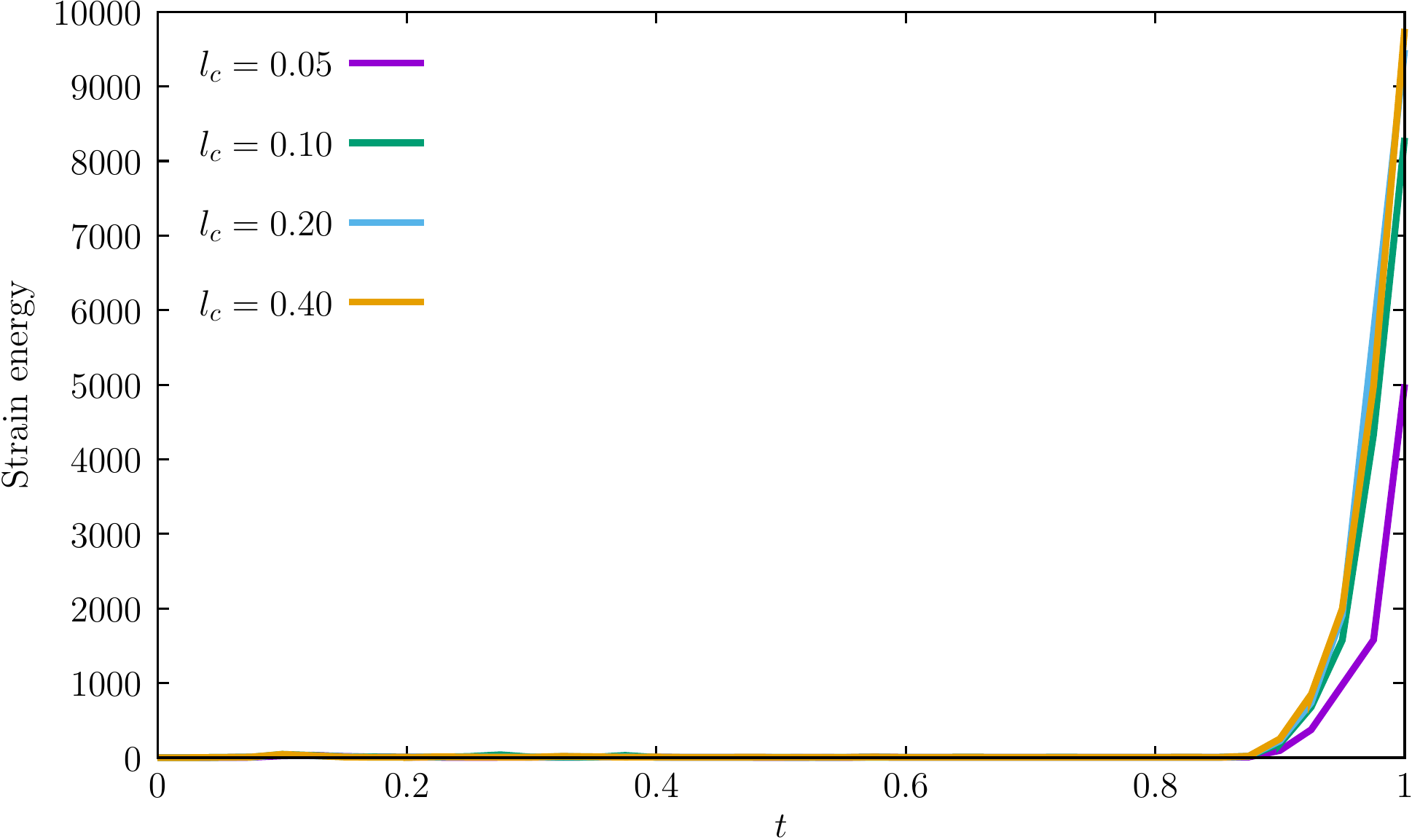}
\par\end{centering}
}\subfloat[]{\begin{centering}
\includegraphics[width=0.50\textwidth]{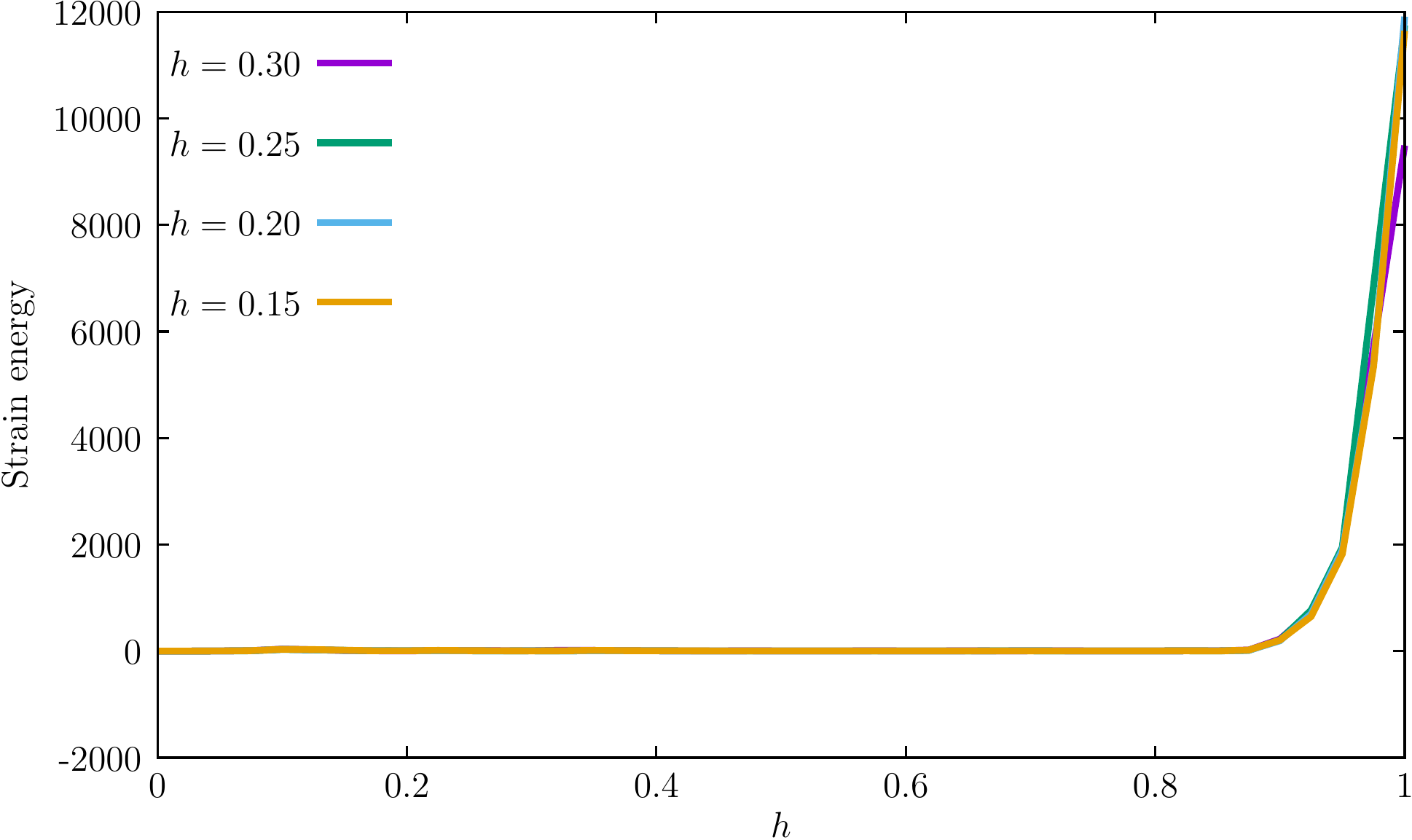}
\par\end{centering}
}
\caption{(a) Effect of $l_{c}$ on the strain energy ($h=0.3$) and
         (b) Effect of $h$ on the strain energy $(l_{c}=0.2$).
         $\kappa=0.6\times10^{6}$ is used.}
\label{fig:Effect-of--3}
\end{figure}

The analysis of gap violation, $v_{\max}=\sup_{\bm{x}\in\Omega}\left[-\min\left(0,g\right)\right]$
as a function of pseudotime $t\in[0,1]$ is especially important for
assessing the robustness of the algorithm with respect to parameters
$l_{c}$ and $\kappa$. For the interval $l_{c}\in[0.05,0.4]$, the
effect of $l_{c}$ is not significant, as can be observed in 
Figure~\ref{fig:Effect-of--1}. Some spikes are noticeable around $t=0.275$
for $l_{c}=0.100$ when the wedge penetrates the cone. Since $\kappa$
is constant, all objects are compressed towards end of the simulation, which the gap violation. In terms of $\kappa$, effects
are the same as in classical geometric-based contact. In terms of
strain energy, higher values of $l_{c}$ result in lower values of
strain energy. This is to be expected, since smaller gradient values
are obtained and the contact force herein is proportional to the product
of the gradient and the penalty parameter. Convergence for the strain
energy as a function of $h$ is presented in Figure~\ref{fig:Effect-of--3}.
It is noticeable that $l_{c}$ has a marked effect near the end of
the compression, since it affects the contact force.
%
%

\subsection{Two-dimensional ironing benchmark}
\begin{figure}
\begin{centering}
\includegraphics[width=0.95\textwidth]{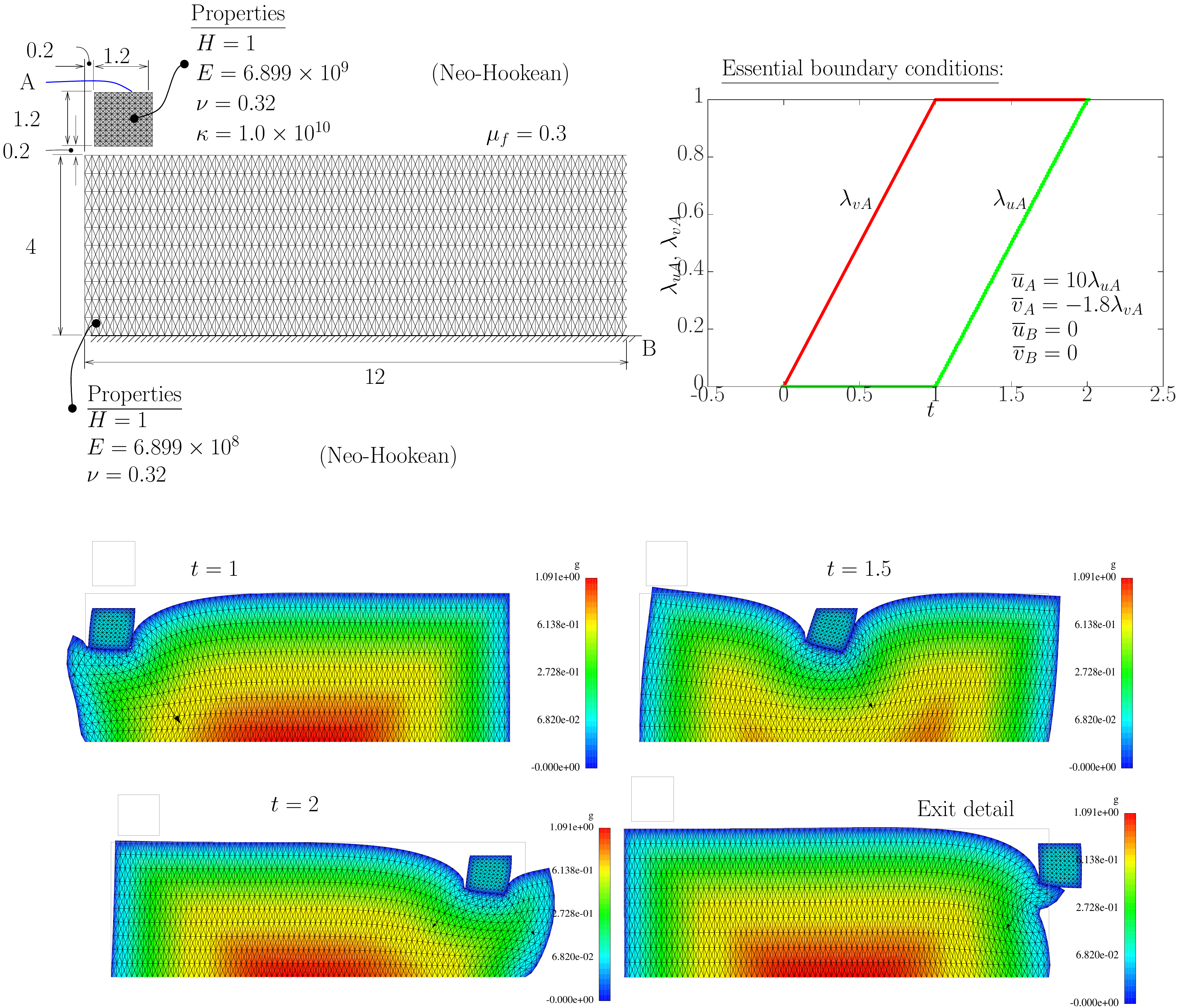}
\par\end{centering}
\caption{\label{fig:Ironing-problem:-relevant}Ironing benchmark in 2D: relevant
data and deformed mesh snapshots.}
\end{figure}

\begin{figure}
\begin{centering}
\subfloat[]{\begin{centering}
\includegraphics[width=0.50\textwidth]{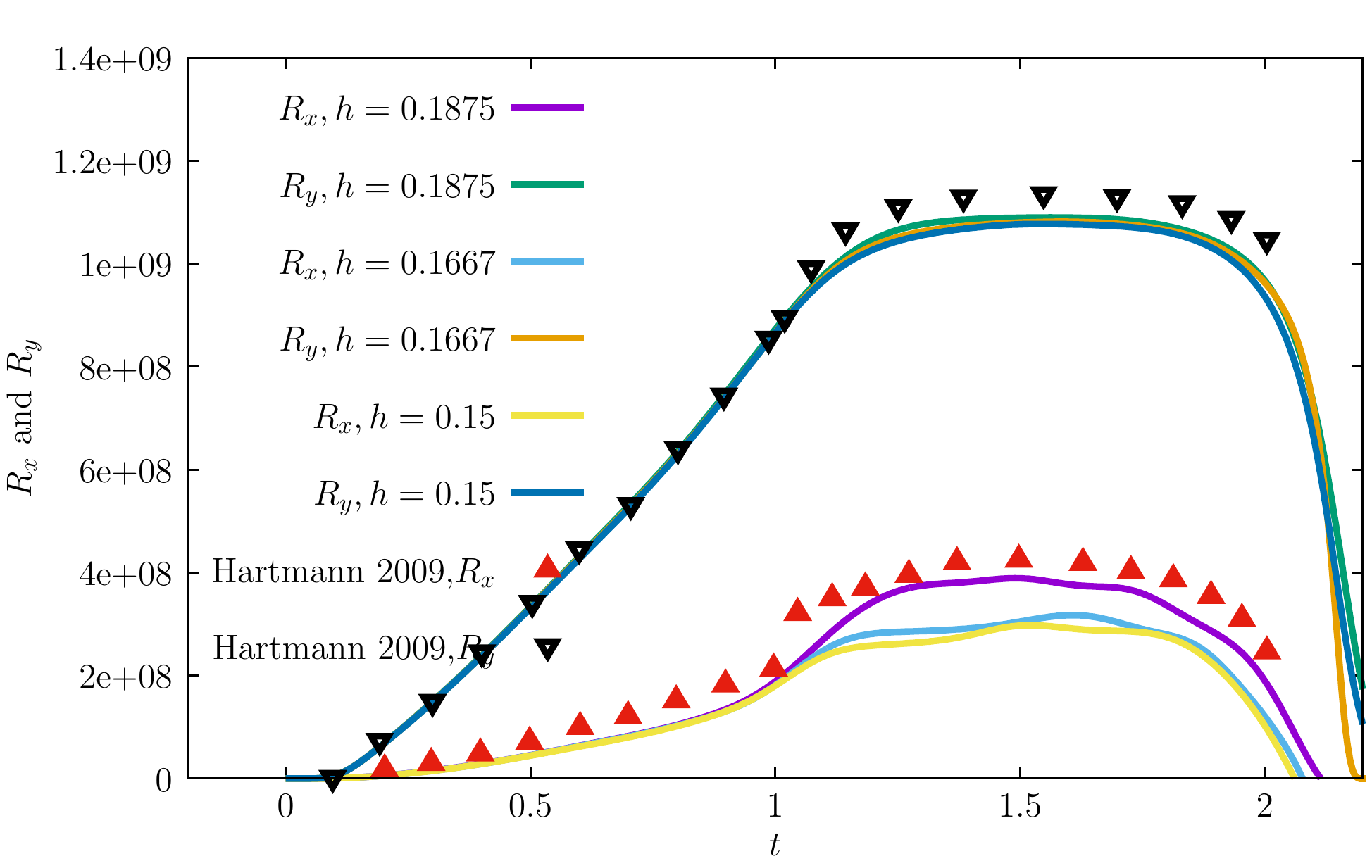}
\par\end{centering}
}\subfloat[]{\selectlanguage{american}%
\begin{centering}
\includegraphics[width=0.50\textwidth]{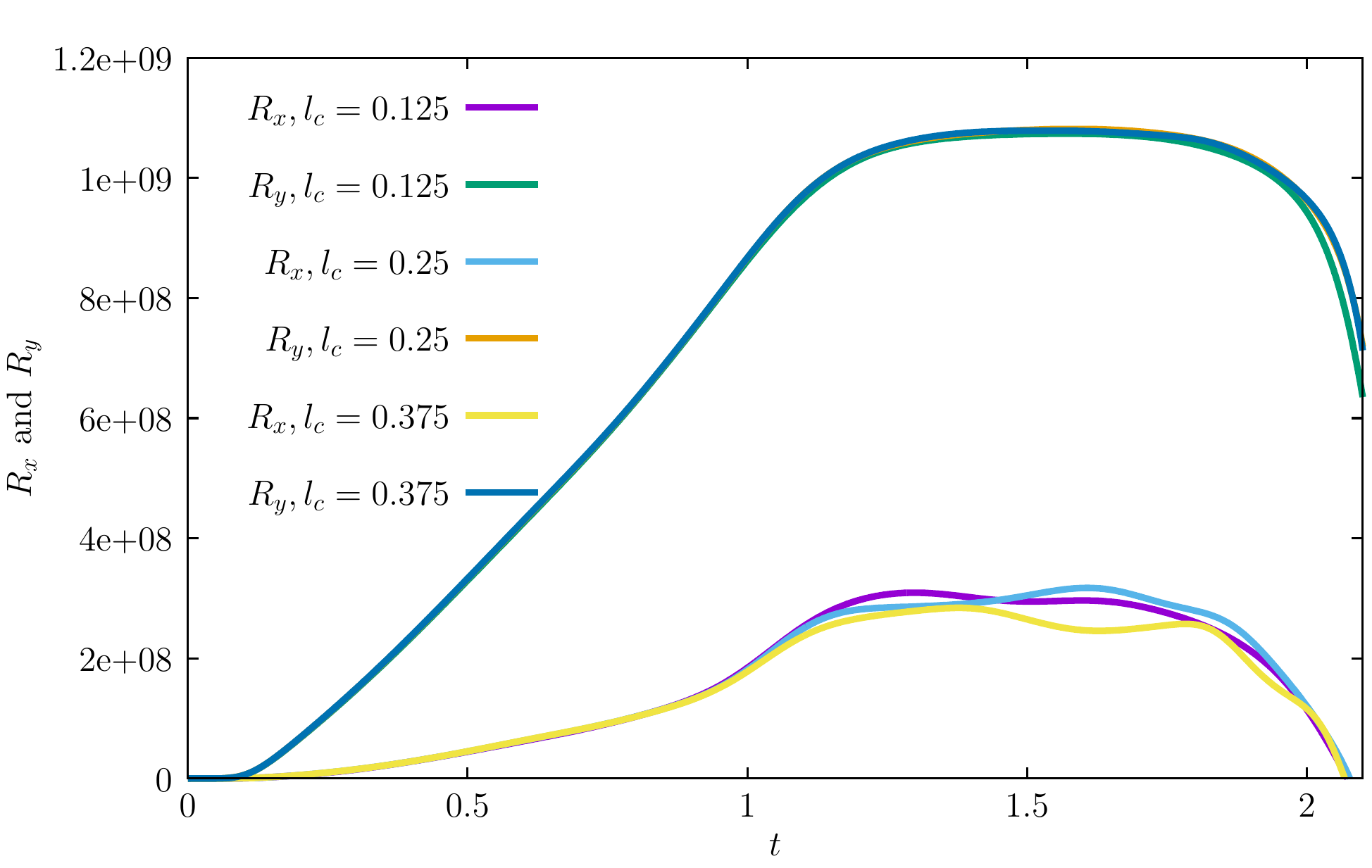}
\par\end{centering}
}
\par\end{centering}
\caption{\label{fig:Ironing-problem:-results}Ironing problem in 2D. Results
for the load in terms of pseudotime are compared to the values reported
in Yang et al.~\cite{yang2005} and Hartmann et al.~\cite{hartmann2009}. $\kappa=0.6\times10^{6}$ is used. (a) Effect of $h$ in the evolution of the horizontal ($R_{x}$) and vertical ($R_{y}$) reactions and (b) Effect of $l_{c}$ in the evolution of the horizontal ($R_{x}$) and
vertical ($R_{y}$) reactions for $h=0.1667$).}
\end{figure}

This problem was proposed by Yang et al.~\cite{yang2005} in
the context of surface-to-surface mortar discretizations. 
Figure~\ref{fig:Ironing-problem:-relevant}
shows the relevant geometric and constitutive data, according to~\cite{yang2005}
and~\cite{hartmann2009}. We compare the present approach with the
results of these two studies in Figure~\ref{fig:Ironing-problem:-results}.
Differences exist in the magnitude of forces, and we infer that 
this is due to the continuum finite element technology. We use finite
strain triangular elements with a compressible neo-Hookean law~\cite{bonet2008}.
The effect of $l_{c}$ is observed in Figure~\ref{fig:Ironing-problem:-results}. 
As before, only a slight effect
is noted in the reaction forces. We use the one-pass version of our
algorithm, where the square indenter has the master nodes and targets
are all elements in the rectangle. Note that, since the cited work
includes friction, we use here a simple model based on regularized
tangential law with a friction coefficient, $\mu_{f}=0.3$. 

\subsection{Three-dimensional ironing benchmark}
\begin{figure}
\selectlanguage{american}%
\begin{centering}
\includegraphics[width=0.85\textwidth]{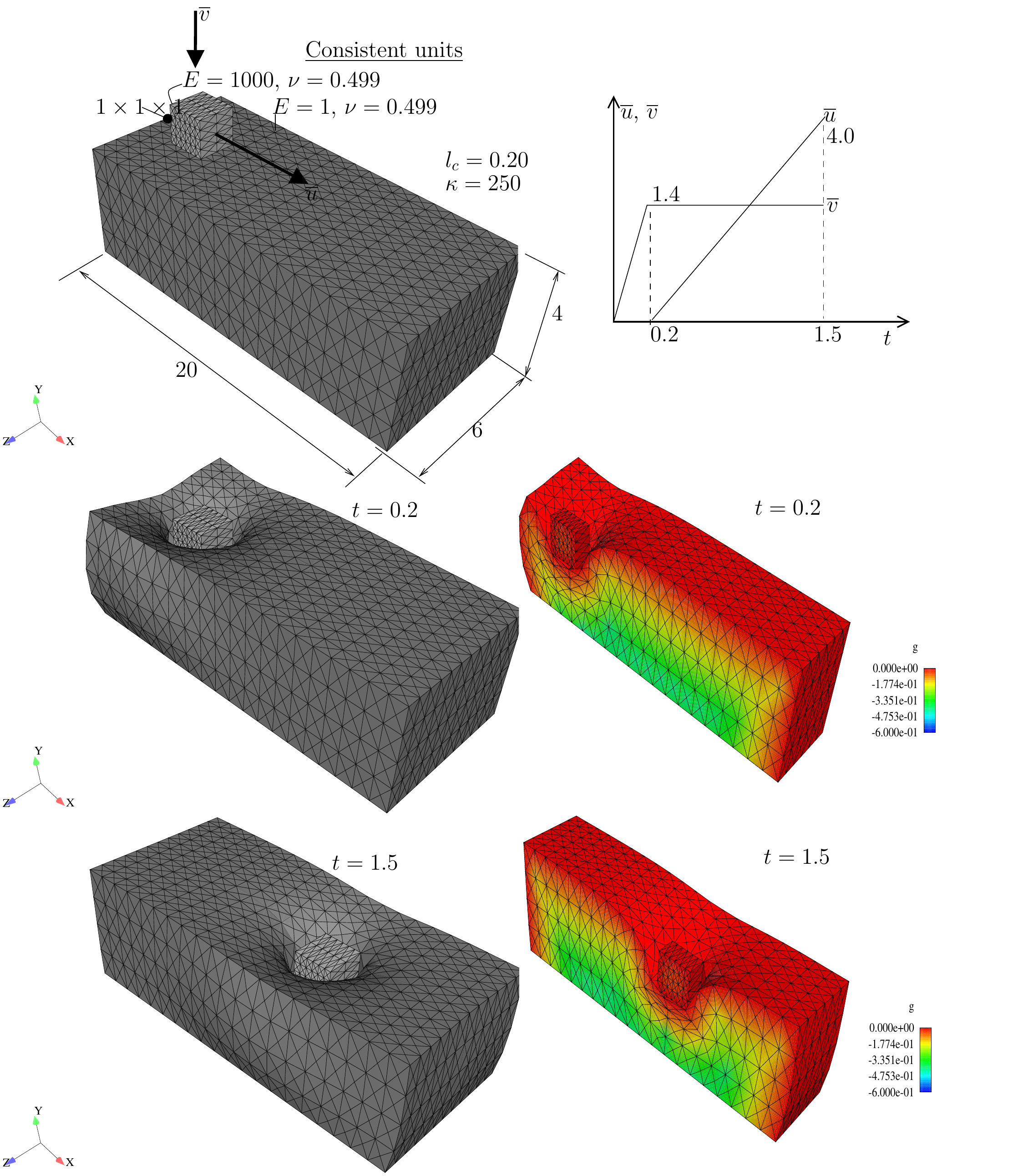}
\par\end{centering}
\selectlanguage{english}%
\caption{\label{fig:I}3D ironing: cube over soft block. Relevant data and
results. }

\selectlanguage{american}%
\selectlanguage{english}%
\end{figure}

\begin{figure}
\begin{centering}
\includegraphics[width=0.7\textwidth]{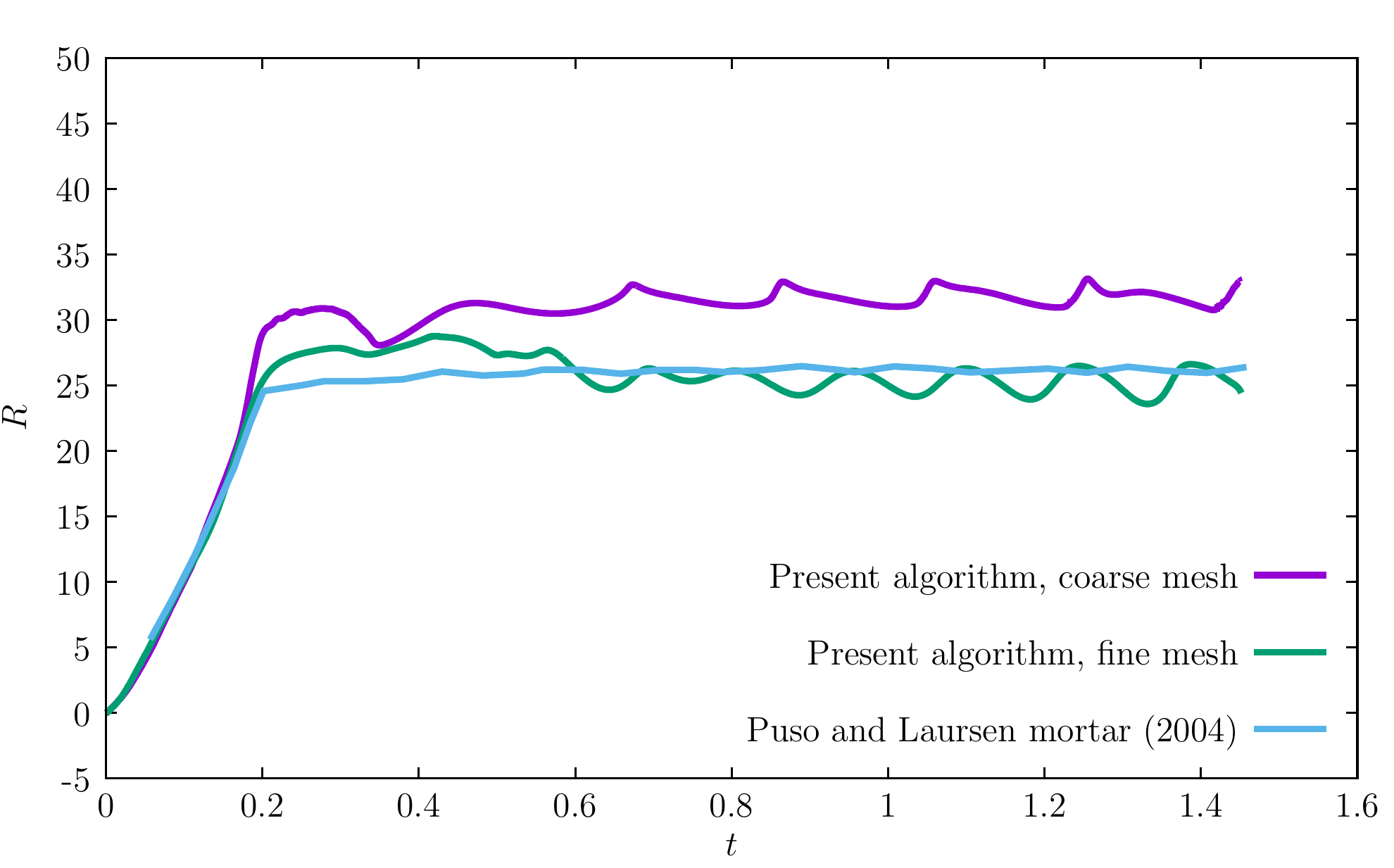}
\par\end{centering}
\caption{\label{fig:3D-ironing:-cube}3D ironing: cube over soft block. Vertical
reactions compared with results in Puso and Laursen~\cite{puso2004b}. }

\end{figure}

We now perform a test of a cubic indenter on a soft block. This
problem was proposed by Puso and Laursen~\cite{puso2004,puso2004b}
to assess a mortar formulation based on averaged normals. The frictionless
version is adopted~\cite{puso2004b}, but we choose the most demanding
case: $\nu=0.499$ and the cubic indenter. Relevant data is presented
in Figure~\ref{fig:I}. The rigid $1\times1\times1$
block is located at $1$ unit from the edge and is first moved down
$1.4$ units. After this, it moves longitudinally $4$ units inside
the soft block. The soft block is analyzed with two distinct meshes:
$4\times6\times20$ divisions and $8\times12\times40$ divisions.
Use is made of one plane of symmetry. A comparison with the vertical
force in~\cite{puso2004b} is performed (see also~\cite{puso2004}
for a clarification concerning the force components). We allow some
interference to avoid locking with tetrahedra. In~\cite{puso2004b},
Puso and Laursen employed mixed hexahedra, which are more flexible
than the crossed-tetrahedra we adopt here. 
Figure~\foreignlanguage{american}{\ref{fig:3D-ironing:-cube}
shows the comparison between the proposed approach and the mortar
method by Puso and Laursen}~\cite{puso2004b}. Oscillations are caused by the normal jumps in gradient of $\phi(\bm{x})$ due to the classical $\mathcal{C}^0$ finite-element discretization (between elements).  Although the oscillations can be observed, the present approach is
simpler than the one in Puso and Laursen.

\selectlanguage{american}%

\clearpage

\section{Conclusions}

We introduced a discretization and gap definition for a contact algorithm
based on the solution of the screened Poisson equation. After a log-transformation, this is equivalent to the solution of a regularized Eikonal equation
and therefore provides a distance to any obstacle or set of obstacles.
This approximate distance function is smooth and is differentiated
to obtain the contact force. This is combined with a Courant--Beltrami
penalty to ensure a differentiable force along the normal direction.
These two features are combined with a step-control algorithm that
ensures a stable target-element identification. The algorithm avoids
most of the geometrical calculations and housekeeping, and is able
to solve problems with nonsmooth geometry. Very robust behavior is
observed and two difficult ironing benchmarks (2D and 3D) are solved
with success. Concerning the selection of the length-scale parameter
$l_{c}$, which produces an exact solution for $l_{c}=0$, we found
that it should be the smallest value that is compatible with the solution
of the screened Poisson equation. Too small of a $l_{c}$ will produce
poor results for $\phi\left(\bm{x}\right)$. Newton-Raphson convergence
was found to be stable, as well as nearly independent of $l_{c}$. In terms of further developments, a $\mathcal{C}^2$ meshless discretization is important to reduce the oscillations caused by normal jumps and  we plan to adopt the cone projection method developed in \cite{areias2015c} for frictional problems.

\label{sec:Conclusions}

\bibliographystyle{unsrt}


\end{document}